\documentclass[12pt,A4]{amsart}
\usepackage{graphicx} 
\usepackage{geometry}
\usepackage{amsmath}
\usepackage{amssymb}
\usepackage{stmaryrd}
\usepackage[hidelinks]{hyperref}

\usepackage{setspace}
\usepackage{graphicx}
\usepackage{float}

\numberwithin{equation}{section}

\theoremstyle{definition}
\newtheorem{theorem}{Theorem}
\newtheorem{proposition}{Proposition}
\newtheorem{lemma}{Lemma}
\newtheorem{sublemma}{sublemma}
\newtheorem{corollary}{Corollary}
\newtheorem{remark}{Remark}

\title{Note on Regularity theory for Nonlinear elliptic equations}
\author{\href{https://zhenye03.github.io}{Zhenye Qian}}
\begin{document}

\begin{abstract}
In this note, we present several seminal developments in the 	\textit{regularity} theory of nonlinear (uniformly) elliptic equations, including the De Giorgi-Nash-Moser theory concerning the Hilbert 19th problem and variational equations, as well as the Krylov-Safonov and Evans-Safonov theories for fully nonlinear equations.
\end{abstract}
\maketitle

\begin{center}
\textit{This note benefited greatly from seminars discussions with my friends. I acknowledge them for their valuable perspectives.}
\end{center}

\section{Introduction}
This note provide a unified exposition of three landmark contributions to the \textit{regularity} theory of second‐order elliptic equations:  
De Giorgi-Nash-Moser theory of divergence‑type equations (1957),  
the Krylov–Safonov Harnack inequality for non‑divergence equations (1979),  
and the Evans–Krylov interior $C^{2,\alpha}$ estimates for concave fully nonlinear equations (1982). All of them make a jump of invariance classes instead of small perturbation.

The general setting is the equation
\[
F(D^{2}u,Du,u,x)=0,
\]
which, under the Schauder theory (\cite{wang2006schauder}) and suitable conditions on the lower‑order terms, can often be reduced to the fully nonlinear model
\[
F(D^{2}u)=0.
\]
We shall always assume that $F$ is \textit{uniformly} elliptic and smooth; the objective is to deduce higher regularity of the solution $u$.  

Historically, the problem originates with Hilbert’s nineteenth problem (1900), which asked whether minimizers of a uniformly convex energy
\[
\mathcal{E}[u]=\int F(\nabla u)\,dx
\]
must be smooth. The Euler–Lagrange equation of such a variational problem is of divergence type,
\[
D_j\left(F_{ij}(\nabla u)D_i u_e\right)=0,
\]
and De Giorgi’s breakthrough \cite{de1957sulla} consisted in combining an energy estimate with the Sobolev inequality in an iterative scheme that yields $C^{1,\alpha}$ regularity, thereby initiating the \textit{bootstrap}.  

For equations that are not in divergence form, the energy estimate is unavailable. Instead, the theory relies on the \textit{maximum principle} and on measure estimates of the contact set, initiated by the works of Calderón–Zygmund (see \cite{muscalu2013classical}) as well as Alexandroff \cite{aleksandrov1958dirichlet}, Bakelman \cite{bakelman1983variational} and Pucci (see \cite{pucci2007maximum}). Building on these ideas, Krylov and Safonov \cite{krylov1979estimate,krylov1981certain} established a \textit{Harnack} inequality for solutions of
\[
a_{ij}(x)D_{ij}u=0,\quad \lambda\delta_{ij}\le a_{ij}\le\Lambda\delta_{ij},
\]
with merely measurable coefficients, which immediately implies $C^{1,\alpha}$ regularity.  

When the equation is both uniformly elliptic and concave (or convex), Evans \cite{evans1980solving,evans1982classical} and Krylov \cite{krylov1983boundedly,krylov1984boundedly} established interior $C^{2,\alpha}$ regularity. Their proof hinges on the observation that second derivatives, say $v$, behave as "almost supersolutions." This property allows one to apply the Krylov–Safonov $L^\epsilon$-estimate to deduce oscillation decay for $v$, which in turn yields the $C^{2,\alpha}$ estimate. Bootstrap arguments then imply that the solution is in fact smooth.

These three theories--De Giorgi–Nash–Moser, Krylov–Safonov, and Evans–Krylov--form the core of modern elliptic regularity theory and have inspired countless extensions to degenerate, nonlocal, and geometric settings.

\medskip

\noindent\textbf{Structure of the notes.}
\begin{enumerate}
    \item[1.] De Giorgi–Nash–Moser theory for divergence‑type equations.
    \item[2.] Krylov–Safonov theory for non‑divergence equations.
    \item[3.] Evans–Krylov theory for concave fully nonlinear equations.
\end{enumerate}

The presentation follows the spirit of several well‑known lecture notes, including those of Caffarelli \cite{caffarelli2010giorgi}, Vasseur \cite{vasseur2016giorgi}, professor and \href{https://sites.math.washington.edu/~yuan/}{Yu Yuan}'s lecture notes, as well as the classical monograph \cite{caffarelli1995fully}. The material was also shaped by a series of inspiring lectures given by Professor Jiakun Liu at the \href{https://bicmr.pku.edu.cn/content/show/17-3335.html}{2024 summer school of BICMR}.

\section{Hilbert 19th problem and De Giorgi-Nash-Moser theory}
\subsection{Introduction to Hilbert 19th problem and calculus of variations}
We begin with Hilbert's 19th problem, proposed at the 1900 ICM \cite{hilbert1900mathematische}. Consider the functional
\[
\mathcal{E}[u]:=\int_\Omega F(\nabla u)\mathrm{d}x,\quad u\in H^1(\Omega)
\]
where $F:\mathbb{R}^n\to\mathbb{R}$ is smooth and \textit{uniformly convex}, and $\Omega\subset\mathbb{R}^n$. The question is: Are local minimizers of $\mathcal{E}$ necessarily smooth?\footnote{Hilbert's original statement: ``\textit{there exist partial differential equations whose integrals are all of necessity analytic functions of the independent variables, that is, in short, equations susceptible of none but analytic solutions}.''} This conjecture was motivated by the Dirichlet energy $E[u]:=\int_\Omega|\nabla u|^2$ (with $F(p)=|p|^2$), whose Euler-Lagrange equation $\Delta u=0$ yields harmonic functions. By Weyl's lemma \cite{han2011elliptic}, any $H^1$ solution is smooth, leading to the conjecture for general $F(\nabla u)$, a more challenging case.

Recall that $u\in H^1$ is a \textit{local minimizer} of $\mathcal{E}$ if,
\[
\mathcal{E}[u]\le \mathcal{E}[u+\varphi],\quad \forall \varphi\in C_c^\infty(\Omega).
\]
The functional $\mathcal{E}$ and $F$ are uniformly convex if, 
\begin{equation}\label{condition: uniformly convex}
    0<\lambda I\le D^2F(p)\le \Lambda I,\quad \forall p\in\mathbb{R}^n,\ \text{for some universal positive $\lambda,\Lambda$}.
\end{equation}

Since $u$ minimizes $\mathcal{E}$, consider a perturbation $u+\epsilon\varphi$ and differentiate:
\[
0=\frac{\mathrm{d}}{\mathrm{d}\epsilon}|_{\epsilon=0}\int_\Omega F(\nabla u+\epsilon\nabla\varphi)=\int_\Omega D_jF(\nabla u)D_j\varphi,\quad \forall \varphi\in C_c^\infty(\Omega).
\]
Thus, $u$ weakly satisfies the Euler-Lagrange equation:
\begin{equation}\label{equation: divergence equation 1}
    D_j(D_jF(\nabla u))=0,\quad \text{in}\ \Omega,
\end{equation}
which is of divergence type. Differentiating \eqref{equation: divergence equation 1} in direction $e\in \mathbb{S}^{n-1}$ gives the linearized equation:
\begin{equation}\label{equation: divergence equation 2}
    D_j(F_{ij}(\nabla u)D_iu_e)=0,
\end{equation}
where $F_{ij}=D_{ij}F,\ u_e=D_eu$. By \eqref{condition: uniformly convex}, both equations are \textit{uniformly elliptic} (throughout the whole note).

To establish smoothness, it suffices to show $u\in C^{1,\alpha}$, initiating a \textit{bootstrap} argument. If $u\in C^{1,\alpha}$, then $F_{ij}(\nabla u)\in C^{0,\alpha}$, and by \textit{Schauder} theory\footnote{We refer readers to \cite{han2011elliptic,gilbarg1977elliptic} for more discussions on classical theories.}, $u \in C^{2,\alpha}$. Iterating yields smoothness. Let $v=u_e$ and $a_{ij}(x)=F_{ij}(\nabla u)$. Then \eqref{equation: divergence equation 2} becomes
\begin{equation}\label{equation: divergence equation 3}
    D_j(a_{ij}(x)D_iv)=0,\quad v\in H^1(\Omega),
\end{equation}
with measurable coefficients $a_{ij}$. Note that $u\in C^{1,\alpha}$ iff $v\in C^\alpha$. Thus, Hilbert's problem reduces to whether solutions of \eqref{equation: divergence equation 3}) are in $C^\alpha$.

Early contributions include \cite{bernstein1904nature,schauder1934lineare,morrey1938solutions}. In 1957, De Giorgi \cite{de1957sulla} solved the problem via an iteration method, which will be presented in this paper. Nash \cite{nash1957parabolic,nash1958continuity} treated the parabolic case, and Moser \cite{moser1960new} provided an alternative method to lift integrability to $L^\infty$. These form the \textbf{De Giorgi-Nash-Moser theory} (DNM theory) for equations like \eqref{equation: divergence equation 3}.

We demonstrate the main theorem of this section.
\begin{theorem}[De Giorgi, 1957 \cite{de1957sulla}]\label{theorem: De Giorgi 1}
Let $v\in H^1(B_1)$ be a solution to \eqref{equation: divergence equation 3}, i.e.,
\[D_j(a_{ij}(x)D_iv)=0,\]
where $a_{ij}(x)$ are uniformly elliptic and measurable. Then
\[\|v\|_{C^\alpha(B_{1/2})}\le C\|v\|_{L^2(B_1)},\]
where $C$ depends on $n,\lambda,\Lambda$.
\end{theorem}

The proof is divided into two parts, say the \textit{iteration} and \textit{oscillation} lemma, and are demonstrated in the following two subsections, respectively, based on \cite{caffarelli2010giorgi,vasseur2016giorgi}.

To set the stage for the regularity theory, we first establish the existence and uniqueness of local minimizers \cite{fernandez2023regularity}.

\begin{proposition}[Existence and uniqueness of local minimizer]
Let $\Omega$ be a bounded Lipschitz domain and suppose
\begin{equation}
    \{w\in H^1(\Omega):\ w=g\ \text{on}\ \partial\Omega\}\ne \emptyset.
\end{equation}
Then there exists a unique local minimizer $u\in H^1(\Omega)$ with $u|_{\partial\Omega}=g$ for $\mathcal{E}$. Moreover, $u$ weakly satisfies \eqref{equation: divergence equation 1}.
\end{proposition}

\begin{proof}
By \eqref{condition: uniformly convex}, without loss of generality, we may assume $F$ attains its minimum at $0$ with $F(0)=0$, so $DF(0)=0$ and
\[\lambda|p|^2\le D^2F(p)\le \Lambda|p|^2,\quad \forall p\in\mathbb{R}^n.\]

For existence, let
\[\mathcal{E}_0:=\inf\{\mathcal{E}[w]:\ w|_{\partial\Omega}=g,\ w\in H^1(\Omega)\},\]
which is finite since $\mathcal{E}_0\le\Lambda \|w\|_{\dot{H}^1(\Omega)}<\infty$. Take a minimizing sequence $\{w_k\}\subset H^1$ with $\mathcal{E}[w_k]\to \mathcal{E}_0$. Then $\|\nabla w_k\|_{L^2(\Omega)}\le \lambda^{-1}\mathcal{E}[w_k]<\infty$, and by Poincaré's inequality, $\{w_k\}$ is bounded in $H^1(\Omega)$. The Rellich-Kondrachov theorem \cite{evans1982classical} and Banach-Alaoglu theorem \cite{brezis2011functional} yield a subsequence $\{w_{k_j}\}$ converging strongly in $L^2(\Omega)$ and weakly in $H^1(\Omega)$ to some $w\in H^1(\Omega)$.

To show $w$ is the minimizer, we prove weak lower semi-continuity of $\mathcal{E}$:
\[\mathcal{E}[w]\le \liminf_{j\to\infty}\mathcal{E}[w_{k_j}]=\mathcal{E}_0.\]
Define $\mathcal{A}(t):=\{v\in H^1(\Omega):\ \mathcal{E}[v]\le t\}$, which is convex. If $v_k\to v$ strongly in $H^1$ with $\{v_k\}\subset \mathcal{A}(t)$, then Fatou's lemma gives
\[\mathcal{E}[v]=\int_\Omega F(\nabla v)\le \liminf_{k\to\infty}\int_\Omega F(\nabla v_k)\le t,\]
so $\mathcal{A}(t)$ is strongly closed. By Mazur's theorem \cite{brezis2011functional}, $\mathcal{A}(t)$ is weakly closed. Now let $t^*:=\liminf_{k\to\infty}\mathcal{E}[w_k]$. For any $\epsilon>0$, there exists a subsequence $\{w_{k_{j,\epsilon}}\}$ with $\mathcal{E}[w_{k_{j,\epsilon}}]\le t^*+\epsilon$. Since $\mathcal{A}(t^*+\epsilon)$ is weakly closed, $\mathcal{E}[w]\le t^*+\epsilon$. Letting $\epsilon\to 0^+$ yields the result.

For uniqueness, suppose $u,v\in H^1(\Omega)$ are two minimizers with $u|_{\partial\Omega}=v|_{\partial\Omega}=g$. If $u\neq v$, there exists $D\subset\Omega$ with $|D|>0$ such that $\nabla u\neq\nabla v$ on $D$. By uniform convexity,
\[\frac{F(\nabla u)+F(\nabla v)}{2}>F\left(\frac{\nabla u+\nabla v}{2}\right)\quad \text{on }D.\]
Integrating over $\Omega$ gives
\[\mathcal{E}_0=\int_\Omega\frac{F(\nabla u)+F(\nabla v)}{2}>\int_\Omega F\left(\frac{\nabla u+\nabla v}{2}\right)\ge \mathcal{E}_0,\]
which leads to the contradiction.
\end{proof}

\subsection{DNM theory I: De Giorgi-Moser iteration}
We now begin the first part of the proof, which establishes an $L^2 \to L^\infty$ estimate. While such an estimate is generally nontrivial, it is expected for solutions of elliptic equations.

\begin{lemma}\label{lemma: DNM 1}
Suppose that
\[\|v^+\|_{L^2(B_1)}<\delta_0\ll 1,\]
where $\delta_0$ is sufficiently small. Then
\[\|v^+\|_{L^\infty(B_1)}\le 1.\]  
\end{lemma}

De Giorgi's iteration combines two key ingredients: the Sobolev inequality and an \textit{energy} inequality, which together yield an iterative scheme by the $L^2$ and $L^\infty$ norms. We present the statements of these lemmas following Caffarelli's Lecture note.

The first lemma applies to general functions:
\begin{sublemma}[Sobolev inequality]
Let $w\in H^1(\mathbb{R}^n)$ have compact support. Then
\[\|w\|_{L^p}\lesssim\|\nabla w\|_{L^2},\]
where $2<p\le \frac{2n}{n-2}$.
\end{sublemma}
\begin{remark}
The lower bound $p>2$ is crucial for the iteration.
\end{remark}
\begin{proof}
We first establish the representation
\begin{equation}\label{equation: representation}
    |w(x_0)|\le \left|\int_{\mathbb{R}^n}\frac{\nabla w(y)\cdot(x_0-y)}{|x_0-y|^n}\mathrm{d}y\right|.
\end{equation}
Without loss of generality, set $x_0=0$. Then
\begin{align*}
    w(0)&=-\frac{\mathrm{d}}{\mathrm{d}t}\frac{1}{|\partial B_1|}\int_{\mathbb{S}^{n-1}}\int_0^\infty w(t\nu)\mathrm{d}t\mathrm{d}\nu\\
    &=-\frac{1}{|\partial B_1|}\int_{\mathbb{S}^{n-1}}\int_0^\infty \nabla w(t\nu)\cdot\nu\mathrm{d}t\mathrm{d}\nu\\
    &=\int_{\mathbb{R}^n}\frac{\nabla w(y)\cdot(0-y)}{|0-y|^n}\mathrm{d}y.
\end{align*}

Define $G(x_0,y)=\frac{x_0-y}{|x_0-y|^n}$. Taking the $L^p$ norm of \eqref{equation: representation} and applying Hölder inequality yields
\begin{align*}
    \int |w(x_0)|^p\mathrm{d}x_0&\le \int \left|\int_{\mathbb{R}^n}\nabla w(y)\cdot G(x_0,y)\mathrm{d}y\right|^p\mathrm{d}x_0\\
    &\le \|\nabla w\|_{L^2}^p\int \|G(x_0,y)\|_{L^2_y}^p\mathrm{d}x_0.
\end{align*}
Since $\|G(x_0,y)\|_{L^2_y}$ is $p$-integrable for $2<p\le \frac{2n}{n-2}$, the result follows.
\end{proof}
\begin{remark}
Replacing the exponent $2$ by a general $q<n$ in the proof yields the Sobolev inequality for $p=\frac{nq}{n-q}$.
\end{remark}

The second lemma concerns subsolutions:
\begin{sublemma}[Energy inequality]
Let nonnegative $v\in H^1$ be a \textit{subsolution} to \eqref{equation: divergence equation 3}, i.e.,
\[D_j(a_{ij}(x)D_iv)\ge0,\quad \text{in}\ B_1\]
Then for any cut-off function $\varphi\in C_c^\infty(B_1)$, we have
\[\int|\nabla (\varphi v)|^2\lesssim C\sup|\nabla\varphi|^2\int(v)^2.\]
\end{sublemma}
\begin{proof}
Multiply the subsolution inequality by $\varphi^2 v\ge 0$ and integrate by parts:
\[\int a_{ij}(2\varphi v D_i\varphi D_jv + \varphi^2 D_iv D_jv) \le 0.\]
Hence
\begin{align*}
    \int a_{ij}D_i(\varphi v)D_j(\varphi v) &= \int a_{ij}\left[(2\varphi v D_i\varphi D_jv + \varphi^2 D_iv D_jv) + D_i\varphi D_j\varphi v^2\right] \\
    &\le \int a_{ij}D_i\varphi D_j\varphi v^2.
\end{align*}
Using the uniform ellipticity $0 < \lambda\delta_{ij} \le a_{ij} \le \Lambda\delta_{ij}$, we obtain
\[\lambda\int |\nabla(\varphi v)|^2 \le \Lambda\sup|\nabla\varphi|^2\int v^2,\]
which completes the proof.
\end{proof}
\begin{remark}
The domain of integration on the right-hand side can be restricted to $B_1 \cap \mathrm{supp}\,\varphi$.
\end{remark}

We now prove Lemma \ref{lemma: DNM 1} using an iterative method.

\begin{proof}[Proof of Lemma \ref{lemma: DNM 1}]
Define cutoff functions $\varphi_k \in C_c^\infty(B_1)$ satisfying
\[\varphi_k:=\begin{cases}
    1,&\text{in}\ B_{1/2+2^{-k}}\\
    0,&\text{outside}\ B_{1/2+2^{-k+1}}
\end{cases},\quad |\nabla\varphi_k|\lesssim 2^k,\ k=1,2,\cdots.\]
Define the \textit{truncated} functions $v_k := (v - (1-2^{-k}))^+$ for $k=1,2,\cdots$. We shall see that $\varphi_k v_k \to 0$ as $k\to\infty$.

By the Sobolev inequality,
\[\left(\int(\varphi_k v_k)^p\right)^{2/p} \lesssim \int|\nabla(\varphi_k v_k)|^2.\]
Since $p > 2$, Hölder inequality yields
\[\int(\varphi_k v_k)^2 \lesssim \left(\int(\varphi_k v_k)^p\right)^{2/p} \cdot |\{\varphi_k v_k > 0\}|^\epsilon,\]
where $\epsilon = 1 - \frac{2}{p} > 0$. Letting $A_k := \int (\varphi_k v_k)^2$, we obtain
\[A_k \lesssim |\{\varphi_k v_k > 0\}|^\epsilon \cdot \int|\nabla(\varphi_k v_k)|^2 =: |\{\varphi_k v_k > 0\}|^\epsilon \cdot E_k.\]

We now estimate $E_k$ using the energy inequality. Note that $\{\varphi_k > 0\} \subset \{\varphi_{k-1} \equiv 1\} \cap \{v_{k-1} > v_k\}$, hence
\[E_k \lesssim \sup|\nabla \varphi_k|^2 \int_{B_1 \cap \mathrm{supp}\,\varphi_k} (v_k)^2 \le 2^{2k} \int (\varphi_{k-1} v_{k-1})^2 = 2^{2k} A_{k-1}.\]

Next, observe that $\{\varphi_k v_k > 0\} \subset \{\varphi_{k-1} v_{k-1} > 2^{-k}\}$. Chebyshev inequality gives
\[|\{\varphi_k v_k > 0\}| \le |\{\varphi_{k-1} v_{k-1} > 2^{-k}\}| \le 2^{2k} \int (\varphi_{k-1} v_{k-1})^2 = 2^{2k} A_{k-1}.\]

Combining these estimates yields
\[A_k \lesssim 2^{2k(1+\epsilon)} A_{k-1}^{1+\epsilon},\quad \epsilon > 0.\]
If $A_1 \le \delta_0$ is sufficiently small, then $A_k \to 0$ as $k \to \infty$, completing the proof.
\end{proof}

\begin{figure}[H]
    \centering
    \includegraphics[width=0.7\linewidth]{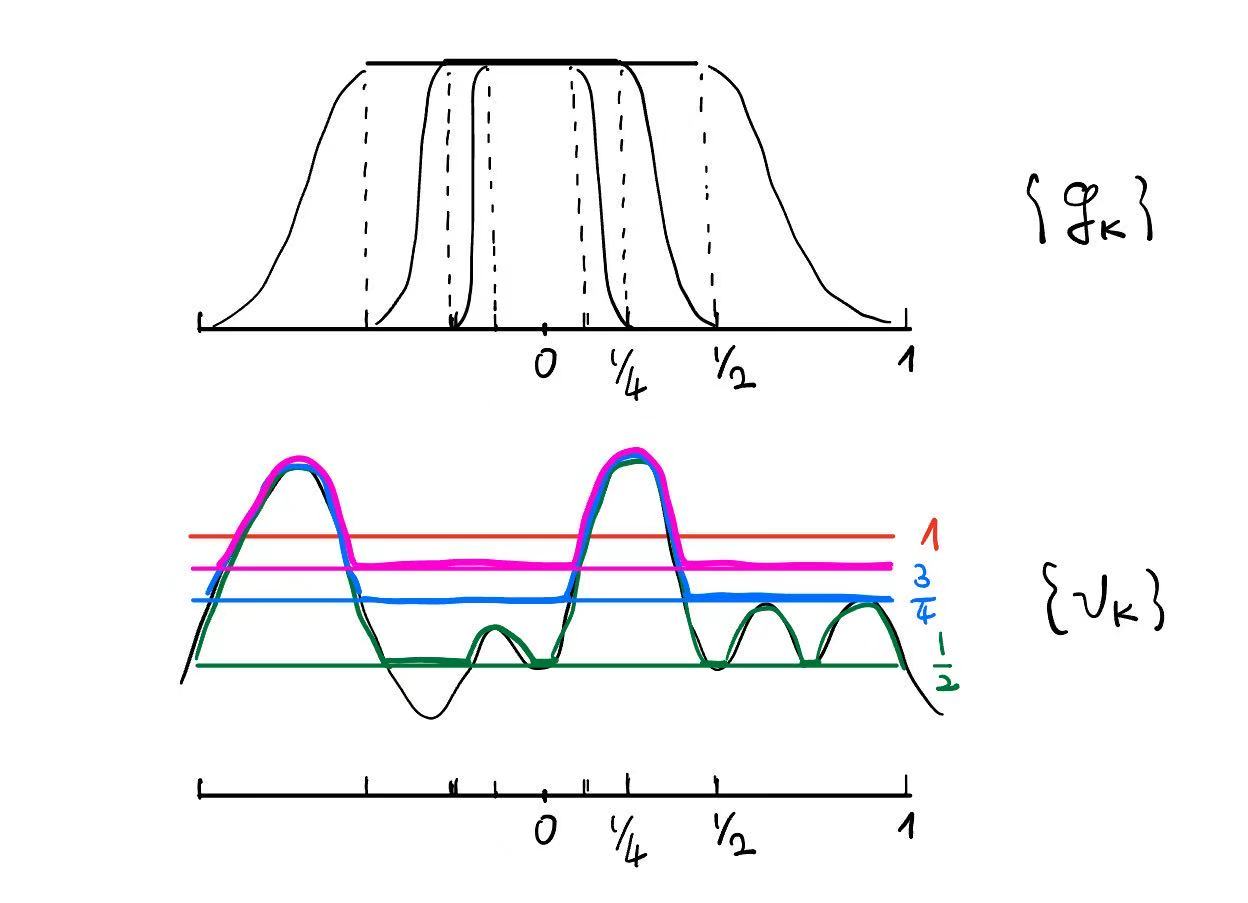}
    \caption{cutoff functions}
\end{figure}

\begin{corollary}
\[\|v\|_{L^\infty(B_{1/2})}\lesssim\|v\|_{L^2(B_1)}\]
\end{corollary}

\begin{remark}
Note that we only use the subsolution $v$.
\end{remark}

\subsection{DNM theory II: Oscillation lemma}
Define the oscillation of a function $v$ over $\Omega$ by
\[
\mathrm{osc}_{\Omega} v := \sup_\Omega v - \inf_\Omega v.
\]
\begin{lemma}[Oscillation decay]\label{lemma: DNM 2}
Let $v$ be a solution to \eqref{equation: divergence equation 3} in $B_1$. Then
\[
\mathrm{osc}_{B_{1/2}} v \le \sigma \, \mathrm{osc}_{B_1} v,
\]
where $\sigma \in (0,1)$ is a universal constant.
\end{lemma}

\begin{remark}
This lemma implies $v \in C^\alpha(B_{1/2})$, completing the proof of Theorem \ref{theorem: De Giorgi 1}.
\end{remark}

Lemma \ref{lemma: DNM 1} ensures that when $\|v\|_{L^2(B_1)}$ is sufficiently small, we have $\mathrm{osc}_{B_1} v \le 2$; in particular, $\mathrm{osc}_{B_{3/4}} v \le 2$. Now suppose $v^+ \equiv 0$ in $B_{3/4}$ except on a set of measure at most $\delta_0/2$, which implies $\|v^+\|_{L^2(B_{3/4})} \le \delta_0/2$. Then Lemma \ref{lemma: DNM 1} yields $\mathrm{osc}_{B_{1/2}} v \le 1$. Note that either $v^+$ or $v^-$ must vanish on at least half of $B_{3/4}$, with loss of generality,
\[\frac{|\{v^+\equiv 0\}\cap B_{3/4}|}{|B_{3/4}|}\ge 1/2.\]
Thus, it remains to establish the vanishing of $v^+$ on a sufficiently large subset. The idea is to iteratively \textit{truncate} and \textit{renormalize} $v^+$ so that the function vanishes on a progressively larger set, ultimately reaching the desired large measure.

\begin{figure}[H]
    \centering
    \includegraphics[width=0.85\linewidth]{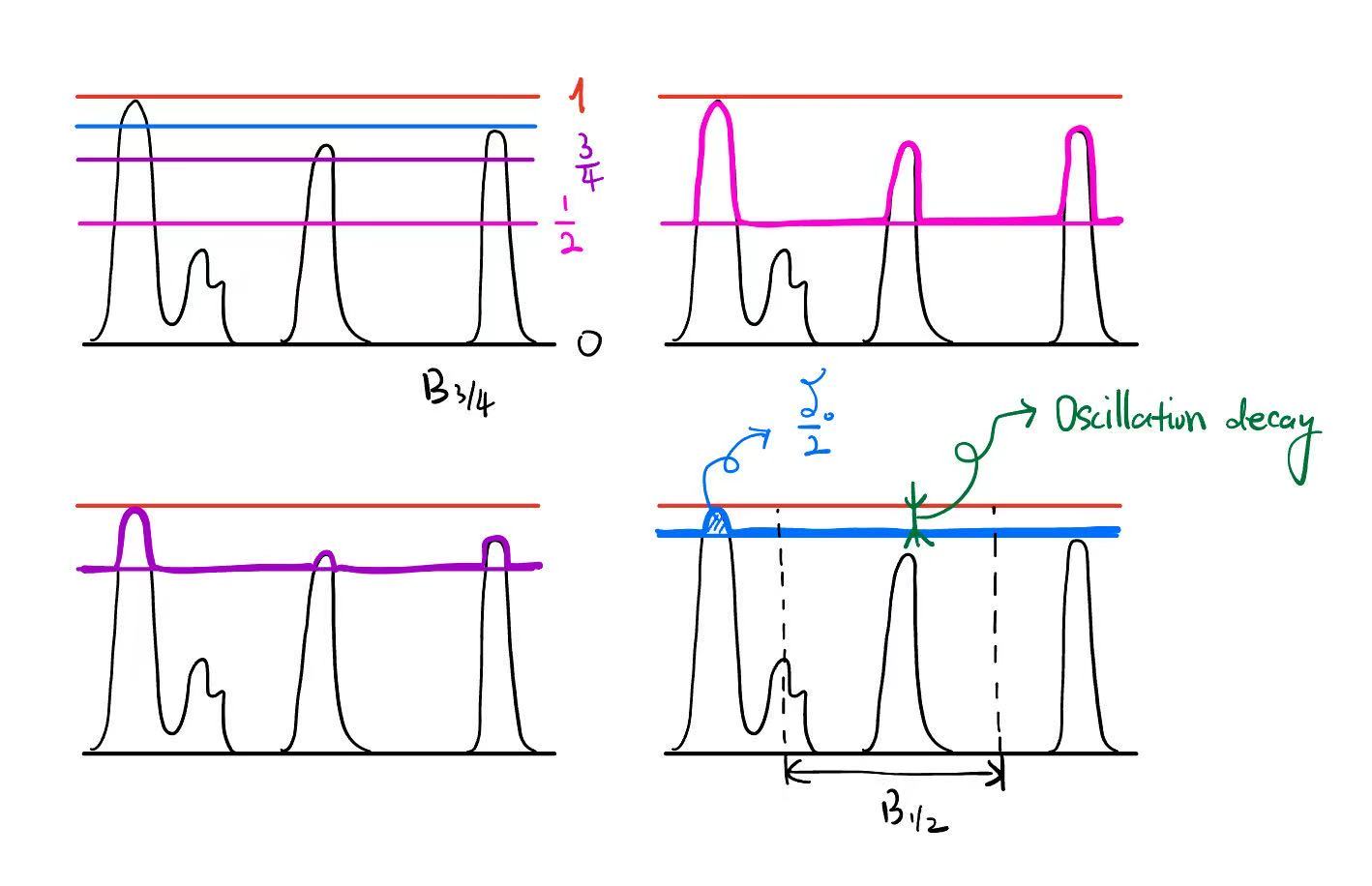}
    \caption{truncate $v^+$}
\end{figure}

First, we shall show that for any $H^1$ function, the set where it remains bounded away from both its supremum and infimum has a quantitatively positive measure, enabling such truncation and renormalization.

\begin{sublemma}[De Giorgi isoperimetric inequality]
Let $w \in H_0^1(B_1)$ satisfy
\[\int_{B_1} |\nabla w|^2 \le 1, \quad 0 \le w \le 1,\quad \text{$w$ touches $1$ at some points.}\]
Define the sets
\[A := \{w = 0\} \cap B_1, \quad C := \{w = 1\} \cap B_1, \quad D := \{0 < w < 1\} \cap B_1.\]
Then
\[|D| \gtrsim |C|^{2 - \frac{2}{n}}.\]
\end{sublemma}

\begin{figure}[H]
    \centering
    \includegraphics[width=0.6\linewidth]{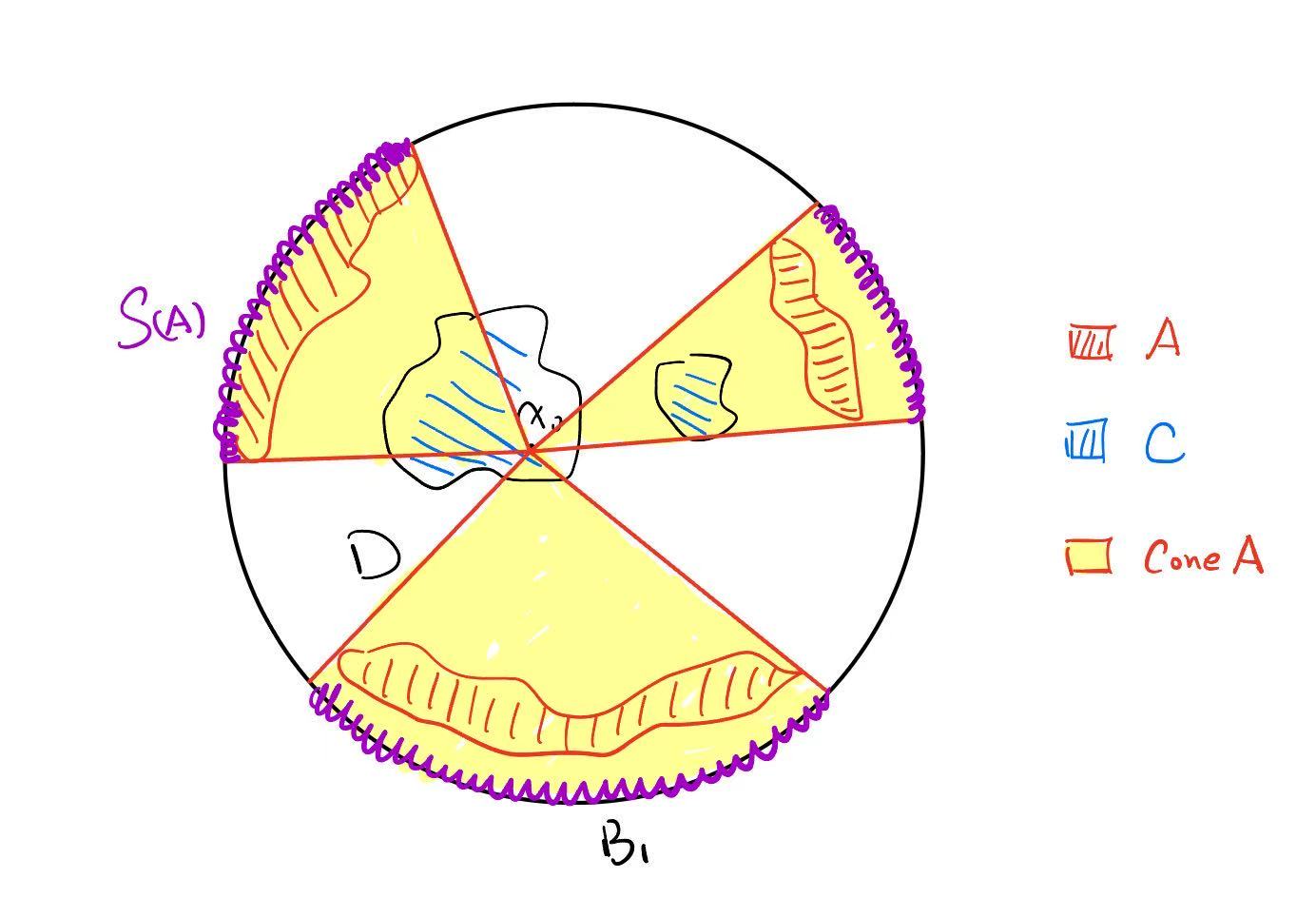}
    \caption{measure estimate}
    \label{measure_1}
\end{figure}

\begin{proof}
Fix $x_0 \in A$ and let $\mathrm{cone}\ A$ be the cone over $A$ with vertex $x_0$ in $B_1$, and let $S(A)$ denote the surface area of $\mathrm{cone}\ A \cap \partial B_1$. See Figure~\ref{measure_1}. Then
\[1 = w(x_0) = \frac{1}{S(A)} \int_{\mathrm{cone}\ A} \frac{\nabla w(y) \cdot (x_0 - y)}{|x_0 - y|^n} \, \mathrm{d}y,\]
or
\[S(A) \le \int_{\mathrm{cone}\ A} \frac{|\nabla w(y)|}{|x_0 - y|^{n-1}} \, \mathrm{d}y.\]
By the \textit{isoperimetric} inequality, $|A|^{(n-1)/n} \lesssim S(A)$. Since $|\nabla w|$ vanishes on $A \cup C$, we have
\[|A|^{(n-1)/n} \lesssim \int_D \frac{|\nabla w(y)|}{|x_0 - y|^{n-1}} \, \mathrm{d}y.\]
Integrating over $x_0 \in C$ and applying Fubini theorem yields
\[|A|^{(n-1)/n} \cdot |C| \lesssim \int_C \int_D \frac{|\nabla w(y)|}{|x_0 - y|^{n-1}} \, \mathrm{d}y \, \mathrm{d}x_0 
= \int_D |\nabla w(y)| \left( \int_C \frac{1}{|x_0 - y|^{n-1}} \, \mathrm{d}x_0 \right) \mathrm{d}y.\]

Note that
\[\int_C \frac{1}{|x_0 - y|^{n-1}} \, \mathrm{d}x_0 \le \frac{1}{|\partial B_1|} \int_{\mathbb{S}^{n-1}} \int_0^{|C|^{1/n}} \frac{1}{r^{n-1}} r^{n-1} \, \mathrm{d}r \, \mathrm{d}\theta = |C|^{1/n},\]
and by Hölder inequality,
\[\int_D |\nabla w| \, \mathrm{d}y \le \|\nabla w\|_{L^2(B_1)} \cdot |D|^{1/2} \le |D|^{1/2}.\]
Thus,
\[|A|^{(n-1)/n} \cdot |C| \lesssim |D|^{1/2} \cdot |C|^{1/n}.\]
Noticing that $|A| \gtrsim |B_1|/2$, we obtain
\[|D| \gtrsim |C|^{2 - \frac{2}{n}},\]
completing the proof.
\end{proof}

\begin{remark}
The proof extends to $W^{1,p}$ functions for $1<p<\infty$.
\end{remark}

Finally, we prove Lemma \ref{lemma: DNM 2}.
\begin{proof}[Proof of Lemma \ref{lemma: DNM 2}]
Assume without loss of generality that $\sup_{B_1}v^+\le 1$. Define the rescaled functions
\[v_k := \frac{v^+ - (1 - 2^{-k})}{2^{-k}}, \quad k = 0,1,2,\cdots\]
and the corresponding sets
\[A_k := \{v_k = 0\} \cap B_{3/4}, \quad C_k := \{v_k = 1\} \cap B_{3/4}, \quad D_k := \{0 < v_k < 1\} \cap B_{3/4}.\]

We prove that for some $k_0$, the function $v_{k_0}$ vanishes in $B_{3/4}$ except on a set of measure at most $\delta_0/2$. Suppose otherwise, and assume $|C_1|\ge \delta>0$. Then by iteration,
\[|A_1|=|A_0|+|D_1|\gtrsim \frac{1}{2}|B_{3/8}|+\delta^{2-\frac{2}{n}},\]
\[|A_k|\gtrsim \frac{1}{2}|B_{3/8}|+k\delta^{2-\frac{2}{n}},\quad k=2,3,\cdots,\]
which yields a contradiction for sufficiently large $k$.
\end{proof}

Finally, we establish a \textit{Liouville-type} result as a direct consequence.
\begin{corollary}\label{corollary: Liouville 1}
Let $v$ be an entire solution of \eqref{equation: divergence equation 3} in $\mathbb{R}^n$, with
\[\|v\|_{L^\infty}\le C.\]
Then $v$ is constant.
\end{corollary}

\begin{proof}[Sketch of proof]
Iterating the oscillation estimate yields
\[
\mathrm{osc}_{B_1} v \le \sigma^k \mathrm{osc}_{B_{2^k}} v \to 0 \quad \text{as } k \to \infty,
\]
completing the proof.
\end{proof}

\begin{remark}
Thus, if the \textit{global} minimizer $u$ satisfies $\|\nabla u\|_{L^\infty}\le C$, then $u$ is a plane.
\end{remark}

\subsection{Density estimate for minimal surfaces}
De Giorgi's iteration scheme proves powerful in establishing decay estimates when two quantities of different homogeneities compete, such as integrability versus differentiability in function theory, or \textit{volume} versus \textit{perimeter} in minimal surfaces and modern geometric measure theory.

As an example, we present an iteration argument in the geometric setting of minimal surfaces. First, recall a well-known density estimate for minimal surfaces \cite{colding2011course,giusti1984minimal}.

\begin{theorem}[Density estimate]
Let $0\in\partial E$ be a minimal hypersurface in $\mathbb{R}^{n+1}$. Then there exists a universal constant $\delta_0>0$ such that
\[|B_1\cap E|\ge \delta_0.\]
\end{theorem}

\begin{proof}
We argue by contradiction. Suppose $|B_1\cap E|<\delta_0$. We show that $B_{1/2}\cap E=\emptyset$.

First, we recall two key tools:
\begin{enumerate}
    \item[A] (\textbf{Isoperimetric inequality}) Since $\partial E$ is area-minimizing in $B_1$, the isoperimetric inequality gives
    \begin{equation}\label{inequality: isoperimetric inequality 1}
        |B_r\cap E|\lesssim P_{B_r}(E)^{\frac{n}{n-1}}\le |2S|^{\frac{n}{n-1}},
    \end{equation}
    where $P_{B_r}(E)$ denotes the perimeter of $E$ in $B_r$, and $S:=S_r$ is the intersection of $E$ with $\partial B_r$.
    \item[B] (\textbf{Excess estimate}) Let $S_\rho$ denote the $\rho$-slice of $(B_R\setminus B_r)\cap E$. Then
    \[|(B_R\setminus B_r)\cap E|=\int_r^R|S_\rho|\mathrm{d}\rho\ge \min_\rho |S_\rho|\cdot (R-r),\]
    so that
    \begin{equation}\label{inequality: energy inequality 1}
        \min_\rho |S_\rho|\le \frac{1}{R-r}|(B_R\setminus B_r)\cap E|.
    \end{equation}
\end{enumerate}

\begin{figure}[H]
    \centering
    \includegraphics[width=0.65\linewidth]{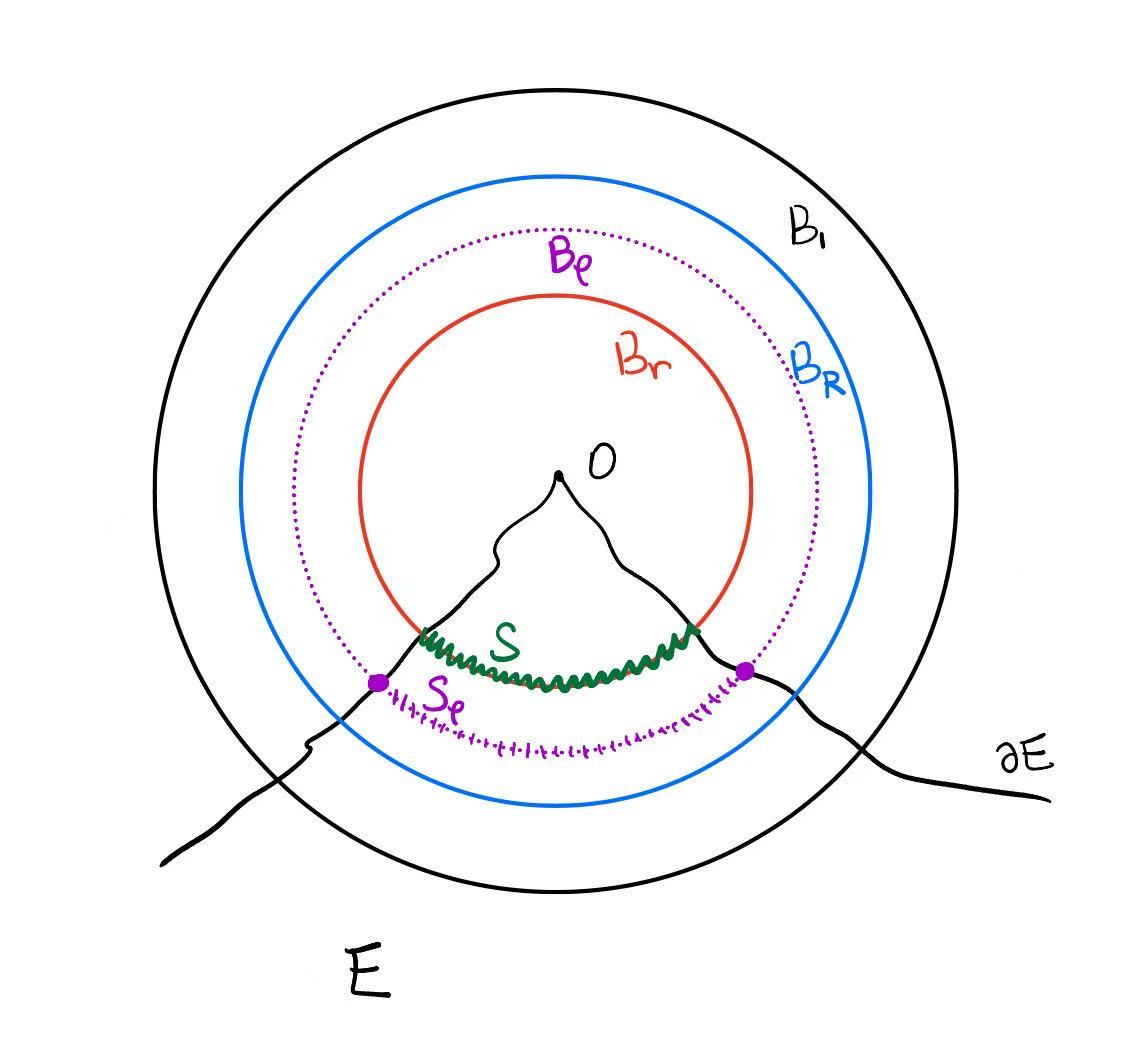}
    \caption{Minimal surface.}
\end{figure}

We now show that $|B_{1/2+2^{-k}}\cap E|\to 0$ as $k\to\infty$. Let $S_k$ be the area-minimizing $\rho$-slice between $B_{1/2+2^{-k-1}}$ and $B_{1/2+2^{-k}}$, and let $\mathcal{C}_k$ be the cone with vertex $0$ and base $S_k$ (see Figure \ref{cutoff surfaces 1}). Applying \eqref{inequality: isoperimetric inequality 1} and \eqref{inequality: energy inequality 1} yields
\[|(B_{1/2+2^{-k}}\setminus B_{1/2+2^{-k-1}})\cap E|\le |\mathcal{C}_{k-1}\cap E|\le |2S_{k-1}|^{\frac{n}{n-1}},\]
and
\[|(B_{1/2+2^{-k}}\setminus B_{1/2+2^{-k-1}})\cap E|\ge |S_k|\cdot(2^{-k}-2^{-k-1})=2^{-k-1}|S_k|.\]

Therefore,
\[|S_k|\le 2^{k+1+\frac{n}{n-1}}|S_{k-1}|^{1+\epsilon},\quad \epsilon=\frac{1}{n-1}.\]
This yields a contradiction for $\delta_0$ sufficiently small, completing the proof.

\begin{figure}[H]
    \centering
    \includegraphics[width=0.65\linewidth]{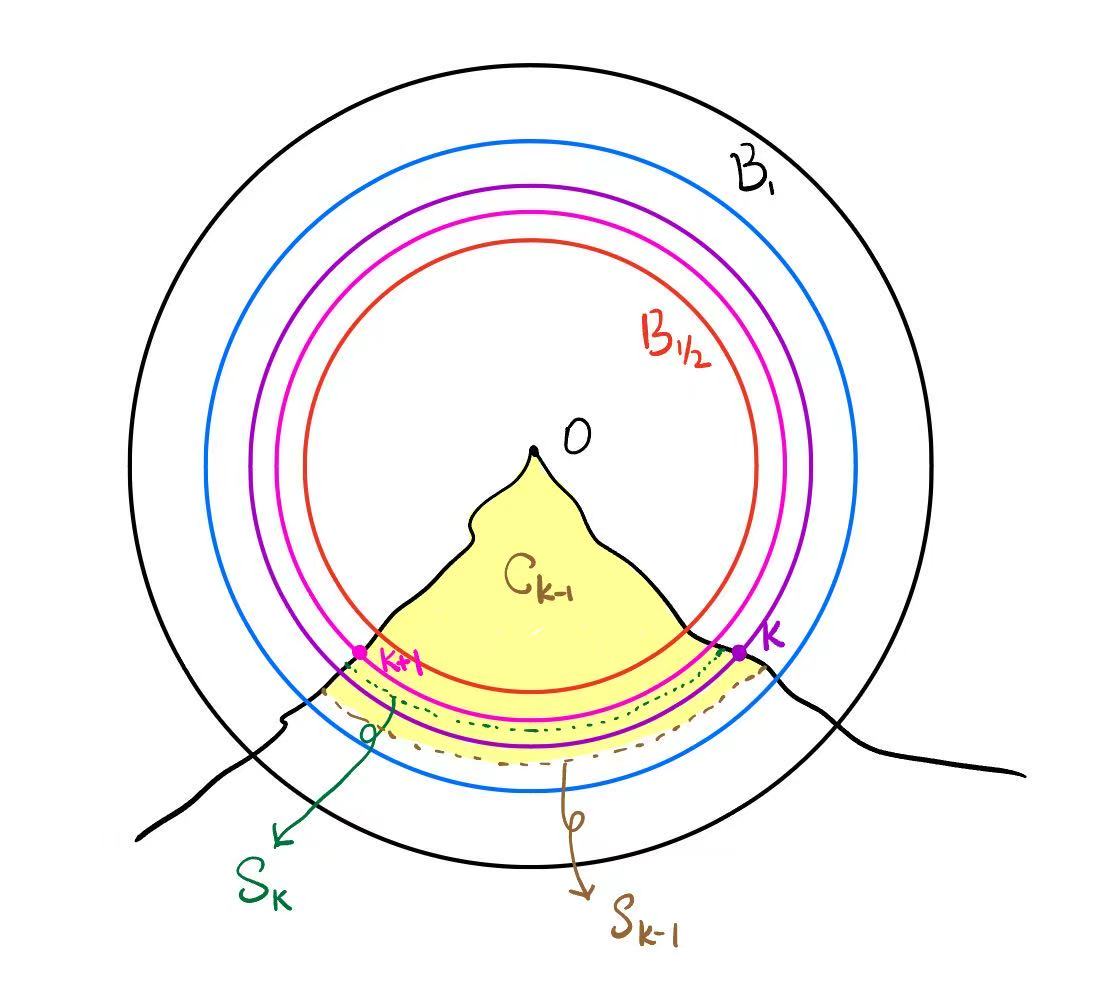}
    \caption{Cut-off surfaces.}
    \label{cutoff surfaces 1}
\end{figure}
\end{proof}

Note that the theorem and the iteration process are analogous to De Giorgi's first step, i.e., Lemma \ref{lemma: DNM 1}. In fact, one may interpret the $L^2$ norm of a function $v$ on a domain $E$ as a volume measure,
\[\int_E (v)^2\sim |E|,\]
while the energy $\|\nabla v\|_{L^2}$ corresponds to a perimeter,
\[\int_E|\nabla v|^2\sim P(E).\]
Thus, Lemma \ref{lemma: DNM 1} can be viewed as a density estimate in the functional setting: 
\[\|v\|_{L^2}\le \delta_0\ll 1 \longleftrightarrow |B_1\cap E|\le \delta_0\ll 1 \Longrightarrow  \|v\|_{L^\infty(B_{1/2})}\le 1\longleftrightarrow B_{1/2}\cap E=\emptyset.\]

This observation similarly applies to general \textit{stationary} varifolds, which are generalized minimal surfaces. For further extensions in geometric measure theory, particularly the Allard regularity theory, we refer the interested reader to \cite{allard1972first} or \cite{fanghua2003geometric}.

\subsection{Moser's alternative approach}
At the end of this section, we briefly outline Moser's alternative approach \cite{moser1960new} to the regularity theory. For a more detailed comparison between Moser’s and De Giorgi’s methods, see \cite{han2011elliptic}. 

Moser’s argument consists of two main parts. The first part achieves the same goal as De Giorgi’s, namely lifting integrability to $L^\infty$, but it replaces the iteration of energy and Sobolev inequalities by a more refined \textit{integration‑by‑parts} technique. It should be noted, however, that the divergence structure of the equation remains essential here; the method does not extend to non‑divergence type equations. The second part is Moser’s \textit{weak} Harnack inequality, which later inspired the Krylov–Safonov theory (compare with Corollary~\ref{corollary: Lepsilon estimate}). Again, this part relies crucially on the divergence form and cannot be directly transferred to the non‑divergence setting.

\paragraph{\textbf{PART I}}
\begin{lemma}[Moser, 1960 \cite{moser1960new}]\label{lemma: moser iteration}
Let $v\ge 0$ be a \textit{subsolution} of \eqref{equation: divergence equation 3} in $B_1$, i.e.,
\[
a_{ij}(x)D_{ij}v\ge 0.
\]
Then for every $p>0$,
\[
\sup_{B_{1/2}} v \le C\|v\|_{L^{p}(B_{1})},
\]
where $C$ depends on $p,n,\lambda,\Lambda$.
\end{lemma}

\begin{proof}
We divide the proof into four steps.
\begin{enumerate}
    \item [\textbf{Step 1}](The case $p\ge 2$) For a nonnegative test function $\varphi\in H^{1}_{0}(B_{1})$ we have
\[
\int a_{ij}D_{i}vD_{j}\varphi\ge 0 .
\]
Choose $\varphi=\eta^{2}v^{p-1}$ with $\eta\in C^{\infty}_{c}(B_{1})$, $\eta\ge0$. Inserting it into the inequality gives
\[
(p-1)\int a_{ij}D_{i}vD_{j}v\cdot v^{p-2}\eta^{2}
\le -2\int a_{ij}D_{i}vD_{j}\eta\cdot v^{p-1}\eta .
\]
Using the uniform ellipticity, we obtain
\[
\int |\nabla v|^{2}v^{p-2}\eta^{2}
\le C \int |\nabla v||\nabla\eta|\,\eta v^{p-1}.
\]
Rewrite the right‑hand side as $\int \left(|\nabla v|v^{\frac{p-2}{2}}\eta\right)
\left(|\nabla\eta|v^{\frac{p}{2}}\right)$,
and apply Young inequality $ab\le\frac{\varepsilon}{2}a^{2}+\frac{1}{2\varepsilon}b^{2}$, we obtain
\[
\int |\nabla v|^{2}v^{p-2}\eta^{2}
\le C_\epsilon\int |\nabla\eta|^{2}v^{p}.\]
Now observe that $|\nabla(v^{p/2})|^{2}= \frac{p^{2}}{4}v^{p-2}|\nabla v|^{2}$; hence
\begin{equation}\label{inequality: moser 1}
    \int |\nabla(v^{p/2})|^{2}\eta^{2}
\le C_1\int |\nabla\eta|^{2}v^{p},
\end{equation}
where $C_1$ is up to $\epsilon,p,\lambda,\Lambda$.
    \item [\textbf{Step 2}](Sobolev embedding)
    Applying the Sobolev inequality to the function $\eta v^{p/2}$ gives
\[
\left(\int  (\eta v^{p/2})^{\frac{2n}{n-2}}\right)^{\!\frac{n-2}{n}}
\le C_{S}\int \left(|\nabla\eta|v^{p/2}+\eta|\nabla(v^{p/2})|\right)^2.
\]
Combining with \eqref{inequality: moser 1} and using Young inequality again, we obtain
\begin{equation}\label{inequality: moser 2}
    \left(\int  v^{\frac{n}{n-2}\cdot p}\cdot\eta^{\frac{2n}{n-2}} \right)^{\frac{n-2}{n}}
\le C_2\int |\nabla\eta|^{2}v^{p},
\end{equation}
where $C_2$ depends on $p,n,\lambda,\Lambda$ for sufficiently small $\epsilon$.
    \item [\textbf{Step 3}](Iteration)
   We define cut‑off functions $\eta_{k}\in C^{\infty}_{c}(B_1)$ such that
\[
\eta_k:=\begin{cases}
    1,&\text{in}\ B_{1/2+2^{-k}},\\
    0,&\text{outside}\ B_{1/2+2^{-k+1}}
\end{cases},\quad k=1,2,\cdots.
\]
Set $p_{k}=p\left(\frac{n}{n-2}\right)^{k}$ (so that $p_{k}\to\infty$). Applying \eqref{inequality: moser 2} with $\eta=\eta_{k}$ and exponent $p_{k}$ yields
\[
\int_{B_{1/2+2^{-k}}} v^{p_k}\lesssim C2^{2k(1+\epsilon)}\left(\int_{B_{1/2+2^{-k+1}}}v^{p_{k-1}}\right)^{1+\epsilon},\quad \epsilon=\frac{2}{n-2}>0,
\]
which leads to the same argument in the proof of Lemma~\ref{lemma: DNM 1}. Thus, we obtain
\begin{equation}\label{inequality: moser 3}
    \sup_{B_{1/2}}v\le C_{0}\,\|v\|_{L^{p}(B_{1})}\quad p\ge2.
\end{equation}
    \item [\textbf{Step 4}](The case $0<p<2$) By \eqref{inequality: moser 3} with $p=2$ and Young inequality, we obtain
    \[\sup_{B_{1/2}}v\lesssim \epsilon\sup_{B_1}v^{1-\frac{p}{2}}+C_\epsilon\left(\int v^p\right)^{1/2}.\]
    Finally, by absorbing lemma, we complete the proof.
\end{enumerate}
\end{proof}

\paragraph{\textbf{PART II}}
\begin{lemma}[Moser weak Harnack inequality, 1960 \cite{moser1960new}]\label{lemma: Moser weak Harnack}
Let $v\ge 0$ be a \textit{supersolution} of \eqref{equation: divergence equation 3} in $B_{1}$, i.e.,
\[
a_{ij}(x)D_{ij}v\le 0.
\]
Then there exist constants $p_{0}>0$ and $C>0$, depending only on $n,\lambda,\Lambda$, such that
\[
\inf_{B_{1/2}} v \ge C\,\|v\|_{L^{p_{0}}(B_{1})}.
\]
\end{lemma}

\begin{proof}[Sketch of proof]
Assume without loss of generality that $v>0$ in $B_{1}$. Apply Lemma~\ref{lemma: moser iteration} to $v^{-1}$, which is nonnegative and satisfies $a_{ij}D_{ij}(v^{-1})\ge0$ (\textit{subsolution}). For any $p>0$,
\begin{equation}\label{inequality: moser 5}
\inf_{B_{1/2}} v \gtrsim\left(\int  v^{-p}\right)^{\!-1/p}.
\end{equation}
Note that the function $w:=-\log v$ is a \textit{subsolution}. By the John-Nirenberg inequality \cite{john1961functions}, there exists $p_{0}=p_{0}(n,\lambda,\Lambda)>0$ such that
\[\int  e^{p_{0}|w-\bar{w} |}\le C,\]
where $\bar{w} =\frac{1}{|B_{1}|}\int w$. Take $p=p_0$ in \eqref{inequality: moser 5}.
\end{proof}

\subsection{Notes}
The regularity theory discussed so far relies on the uniform ellipticity condition \eqref{condition: uniformly convex}. A natural question is whether this condition can be relaxed to mere \textit{strict} ellipticity (or \textit{strict} convexity of $F$). Under suitable boundary conditions, minimizers of $\mathcal{E}$ are then Lipschitz, i.e., $\|\nabla u\|_{L^\infty}<\infty$. Whether such minimizers are necessarily $C^1$ remained open for decades. De Silva and Savin \cite{de2010minimizers} answered this positively in dimension $n=2$, by assume that the set of degenerate points is at most finite. For $n\ge 4$, Mooney \cite{mooney2016some,mooney2020minimizers} constructed counterexamples. The case $n=3$ remains open; see \cite{mooney2022hilbert} for further discussion.

The ideas behind De Giorgi's iteration have been extended to many other settings. A prominent variant is the "$\epsilon$-\textit{regularity}" theory: there exists a universal $\epsilon>0$ such that if the energy of a solution is smaller than $\epsilon$, then the solution is bounded. In geometric measure theory, this philosophy appears in Allard's regularity theorem for minimal surfaces, where the relevant energy—called the \textit{excess}—measures the deviation of an approximate tangent plane of the rectifiable set from a fixed direction. On the other hand, Moser's approach was applied by Schoen–Simon–Yau \cite{schoen1975curvature} (see also Theorem 2.21 of \cite{colding2011course}) to obtain $L^p$ estimates for minimal surfaces, replacing the $L^2$ estimates (stability) given by Simons’ inequality \cite{simons1968minimal} (Lemma 2.1 of \cite{colding2011course}). This led to a solution of the stable Bernstein problem under a volume growth assumption for $n\le 5$, with a dimensional gap for the exponent $p$. Recently, \cite{bellettini2025extensions} used De Giorgi’s iteration to extend the result to $n=6$.

The method has since been adapted to various degenerate and nonlocal settings, such as the $p$-Laplacian, degenerate parabolic equations, and fractional diffusion, and has found applications in diverse fields including fluid dynamics \cite{caffarelli2010giorgi} and geometric analysis; see for instance \cite{vasseur2016giorgi, brigati2025introduction, violinimodern} for further discussions.

\section{Fully nonlinear equations and Krylov-Safonov theory}
\subsection{Introduction to fully nonlinear equations}
In the next two sections, we study the regularity theory of \textit{fully nonlinear} equations,
\begin{equation}\label{equation: fully nonlinear equation 1}
    F(D^2u)=0.
\end{equation}
Key results were established by Krylov-Safonov \cite{krylov1979estimate,krylov1981certain} for $C^{1,\alpha}$ regularity, and by Evans-Krylov \cite{evans1982classical,krylov1983boundedly,krylov1984boundedly} for $C^{2,\alpha}$ regularity.

Throughout, we assume $F$ is \textit{uniformly elliptic}: there exist constants $\lambda,\Lambda>0$ such that
\[\lambda\|N^+\|-\Lambda\|N^-\|\le F(M+N)-F(M)\le \Lambda\|N^+\|-\lambda\|N^-\|,\]
which locally is equivalent to 
\[\lambda\delta_{ij}\le F_{ij}\le \Lambda\delta_{ij},\quad F_{ij}:=\frac{\partial F}{\partial u_{ij}}.\]
Assuming $F(0)=0$ without loss of generality, we have
\[\lambda\|M^+\|-\Lambda\|M^-\|\le F(M)\le \Lambda\|M^+\|-\lambda\|M^-\|.\]
Thus, the uniform ellipticity of $F$ implies that for any “solution” $u$ of \eqref{equation: fully nonlinear equation 1}, its Hessian satisfies $\|(D^2u)^+\|\sim\|(D^2u)^-\|$. Consequently, there exist pointwise-defined coefficients $\{a_{ij}(x)\}$ (merely measurable and bounded) such that
\begin{equation}\label{equation: fully nonlinear equation 2}
    a_{ij}(x)D_{ij}u=0,\quad \lambda\delta_{ij}\le a_{ij}\le \Lambda\delta_{ij},
\end{equation}
which is equivalent to the implicit form \eqref{equation: fully nonlinear equation 1}.

However, unlike the divergence‑type equations studied earlier, which admit a natural weak formulation via integration by parts.\footnote{Hence they are often called “integral” or "global" equations, the fully nonlinear equation \eqref{equation: fully nonlinear equation 2}) is inherently local and does not possess a straightforward weak solution concept.} Crandall–Lions \cite{crandall1983viscosity} and Evans \cite{evans1978convergence,evans1980solving} introduced the concept of viscosity solutions for fully nonlinear equations, which serves as a natural replacement for the weak solution notion used in divergence‑type equations. We say that $u$ is a \emph{viscosity supersolution} to \eqref{equation: fully nonlinear equation 2} in a domain $\Omega$ if for any test function $\varphi\in C^2(\Omega)$ and any point $x_0\in\Omega$ where $\varphi$ touches $u$ from below (i.e., $\varphi(x_0)=u(x_0)$ and $\varphi(x)\le u(x)$ for all $x$ in a neighborhood of $x_0$), we have
\[a_{ij}(x_0) D_{ij}\varphi(x_0)\le 0.\]
Analogously, $u$ is a \emph{viscosity subsolution} if for any $\varphi\in C^2(\Omega)$ touching $u$ from above at $x_0$ (i.e., $\varphi(x_0)=u(x_0)$ and $\varphi(x)\ge u(x)$ near $x_0$), we have
\[a_{ij}(x_0) D_{ij}\varphi(x_0)\ge 0.\]
Finally, $u$ is a \emph{viscosity solution} if it is both a viscosity supersolution and a viscosity subsolution. We refer readers to \cite{caffarelli1995fully,fernandez2023regularity} for more details.

The question is whether $u$ enjoys the same $C^{1,\alpha}$ regularity as in De Giorgi's result, and further whether it is smooth. Indeed, differentiating \eqref{equation: fully nonlinear equation 1} in a direction $e\in \mathbb{S}^{n-1}$ yields
\begin{equation}\label{equation: fully nonlinear equation 3}
    F_{ij}(D^2u)D_{ij}u_e=0.
\end{equation}
Krylov-Safonov proved that $u_e\in C^{1,\alpha}$; however, this alone cannot initiate a bootstrap argument because the coefficients $F_{ij}$ depend on second‑order derivatives. Evans-Krylov established higher $C^{2,\alpha}$ regularity under the additional assumption that $F$ is concave, which does start the bootstrap and leads to smooth solutions. Both theories will be presented in this paper.

For simplicity, set $F_{ij}=a_{ij}(x)$, where $a_{ij}$ are merely measurable and bounded, and let $v=u_e$. We prove that $v$, which solves
\begin{equation}\label{equation: fully nonlinear equation 4}
    a_{ij}(x)D_{ij}v=0,\quad \lambda\delta_{ij}\le a_{ij}(x)\le \Lambda\delta_{ij},
\end{equation}
belongs to $C^{\alpha}$. More generally, we state the Krylov-Safonov \textit{Harnack} inequality (\textbf{Krylov-Safonov theory}), whose proof will be divided into two parts.

\begin{theorem}[Krylov-Safonov, 1979 \cite{krylov1979estimate,krylov1981certain}]\label{theorem:KS-Harnack}
Let $v$ be a nonnegative solution to \eqref{equation: fully nonlinear equation 4} in $B_1$, where $a_{ij}$ are \textit{uniformly elliptic}, measurable and bounded. Then
\[
\sup_{B_{1/2}} v \le C \inf_{B_{1/2}} v,
\]
where $C$ depends only on $n,\lambda,\Lambda$.
\end{theorem}

Intuitively, one may view the viscosity solution $v$ as a kind of \textit{harmonic} function.

Observe that the $L^\epsilon$-norm of a superharmonic function is controlled by its value at a point; hence, if we fix the value of $v$ at some point, say $v(0)=1$ without loss of generality, then $v$ belongs to $L^\epsilon$. See Lemma~\ref{lemma: Moser weak Harnack} and Lemma~\ref{lemma: Krylov-Safonov weak Harnack inequality} below.

Conversely, the value of a subharmonic function is controlled by the nearby $L^\epsilon$-norm, so a fixed $L^\epsilon$-norm may govern the value at the center. See Lemma~\ref{lemma: DNM 1}, Lemma~\ref{lemma: moser iteration} and Lemma~\ref{lemma: local maximum principle} below.

Combining these two facts, one can propagate control on $v$ from a larger domain to a smaller one, which yields an estimate on the oscillation. These correspond to the two parts of the proof.

\begin{figure}[H]
    \centering
    \includegraphics[width=0.45\linewidth]{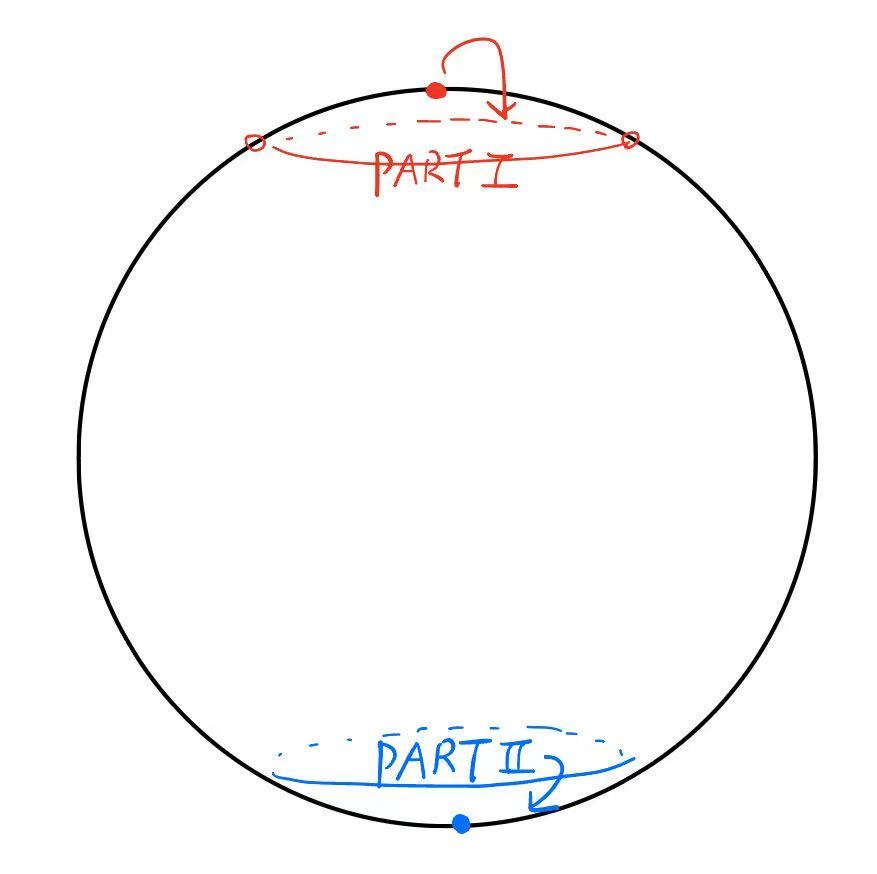}
    \caption{mean valur inequalities}
\end{figure}

\subsection{Two main tools}
Compared with the divergence‑type equation \eqref{equation: divergence equation 3}, the equation \eqref{equation: fully nonlinear equation 4} lacks a crucial ingredient: an energy inequality. This absence is not merely a technical restriction\footnote{Indeed, the energy $|\nabla u|$ is naturally controlled by the functional $\mathcal{E}$. Also, see Lemma~\ref{lemma: moser iteration}.}. Consequently, the iteration argument used for divergence‑type equations fails in the nondivergence case, and one must instead exploit intrinsic structures of elliptic equations, such as the \textit{maximum principle} or comparison principle.

\begin{lemma}[Alexandroff–Bakelman–Pucci, 1960s \cite{aleksandrov1958dirichlet,bakelman1983variational,pucci2007maximum}]\label{lemma: ABP estimate}
Let $v$ satisfy
\[a_{ij}(x)D_{ij}v=f \quad \text{in } B_1,\]
with $v\ge 0$ on $\partial B_1$. Then
\[\sup (v^-)^n \lesssim \int_{\Gamma^+} |f^+|^n,\]
where $\Gamma^+$ denotes the \textit{contact set} of $v$.
\end{lemma}

\begin{remark}
The contact set $\Gamma^+$ is defined as
\[\Gamma^+:=\{y\in B_1:\ p(y)\cdot(x-y)\le v(x)-v(y),\ \text{for some plane $p(y)$ through $y$}\}.\]
The primary utility of the ABP estimate lies in obtaining precise bounds for the \textit{measure} of the contact set.
\end{remark}

\begin{figure}[H]
    \centering
    \includegraphics[width=0.6\linewidth]{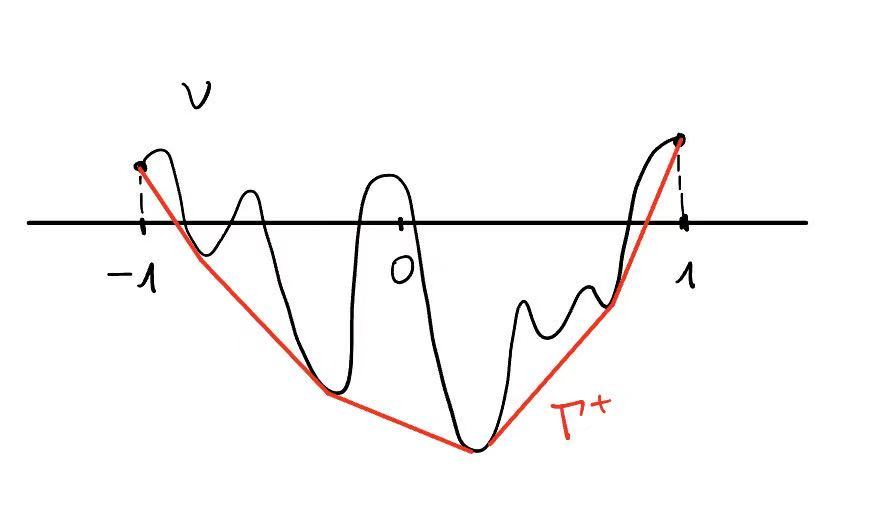}
    \caption{Contact set.}
\end{figure}

A detailed proof can be found in standard PDE textbooks \cite{gilbarg1977elliptic,han2011elliptic}; we provide a brief sketch.
\begin{proof}[Sketch of proof]
First, by a change of variables we have
\[\int_{Dv(\Gamma^+)} g \le \int_{\Gamma^+} g(Dv)|\det D^2 v|, \quad \forall \text{nonnegative } g\in L^{1}_{\mathrm{loc}}.\]
Then, by moving planes to contact points, one can show that the ball $B_{\tilde{M}}(0)$ is contained in $Dv(\Gamma^+)$, where $\tilde{M}\sim \sup v^-$. Choosing an appropriate $g$ completes the proof.
\end{proof}

Second, we state a Calderón–Zygmund type lemma.

\begin{lemma}[A–B lemma]\label{lemma: A-B lemma}
Let $A\subset B\subset \mathbb{Q}_1$, where $\mathbb{Q}_1$ is the unit cube. Assume
\begin{enumerate}
    \item[(i)] $|A|\le \delta$;
    \item[(ii)] For any sub‑cube $Q\subset\mathbb{Q}_1$, if $\frac{|A\cap Q|}{|Q|}>\delta$, then the dilated cube $\tilde{Q}$ (for example, obtained by expanding each edge by a factor of $3$) is contained in $B$.
\end{enumerate}
Then $|A|\le \delta|B|$.
\end{lemma}

\begin{figure}[H]
    \centering
    \includegraphics[width=0.67\linewidth]{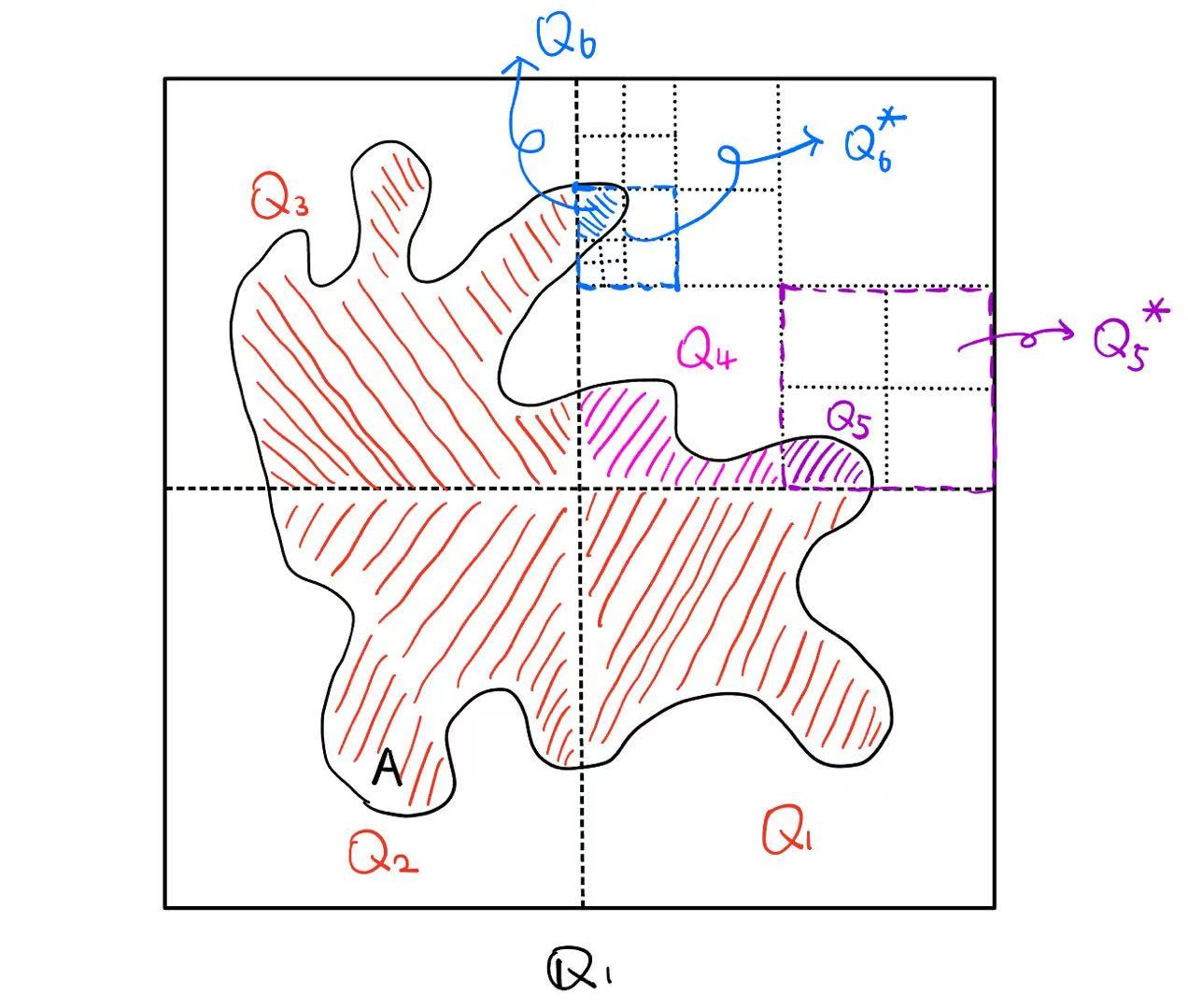}
    \caption{Calderón–Zygmund decomposition.}
    \label{C-Z decomposition}
\end{figure}

\begin{proof}
We perform a Calderón–Zygmund decomposition \cite{muscalu2013classical} of $A$ inside $\mathbb{Q}_1$ (see Figure~\ref{C-Z decomposition}). Let $\{Q_j\}$ be the collection of non‑overlapping dyadic sub‑cubes such that
\[
\frac{|A\cap Q_j|}{|Q_j|}>\delta,
\]
and let $Q_j^*$ denote their predecessors (the last cubes in the decomposition that were not selected). By construction,
\[
\frac{|A\cap Q_j^*|}{|Q_j^*|}\le\delta.
\]
Condition (ii) implies $\tilde{Q}_j\subset B$, hence $Q_j^*\subset B$. After removing repetitions and reindexing, we obtain a family of disjoint cubes $\{Q_{j'}^*\}$ with $\bigcup Q_{j'}^*\subset B$. Consequently,
\[
|A|\le \left|\bigcup Q_j\right|
      \le \left|\bigcup (Q_{j'}^*\cap A)\right|
      \le \sum |Q_{j'}^*\cap A|
      \le \delta\sum |Q_{j'}^*|
      \le \delta|B|,
\]
where the last inequality uses the disjointness of the $Q_{j'}^*$.
\end{proof}

\subsection{Krylov-Safonov weak Harnack inequality}
We now present the first part of the proof of Theorem~\ref{theorem:KS-Harnack}.

\begin{lemma}[Krylov-Safonov weak Harnack inequality]\label{lemma: Krylov-Safonov weak Harnack inequality}
Let $v$ be a supersolution to \eqref{equation: fully nonlinear equation 4} with $v(0)=1$. Then $v\in L^\epsilon_{weak}(B_{1/2})$, i.e.,
\[|\{v>\lambda\}\cap B_{1/2}|\le \lambda^{-\epsilon},\quad \forall \lambda>0.\]
\end{lemma}
\begin{remark}
Recall \textit{Chebyshev} inequality.
\end{remark}

The proof proceeds in two steps.

First, we derive a contact estimate via the ABP method. To present a baby version, We claim that
\[|\{v<2\}\cap B_1|\ge \theta>0,\]
for some universal $\theta>0$. Define
\[\tilde{v}:=v+2(|x|^2-1).\]
Then $\tilde{v}(0)=-1$ and $\tilde{v}|_{\partial B_1}\ge 0$. A direct calculation gives
\[a_{ij}D_{ij}\tilde{v}\le C,\]
for some constant $C>0$.

\begin{figure}[H]
    \centering
    \includegraphics[width=0.7\linewidth]{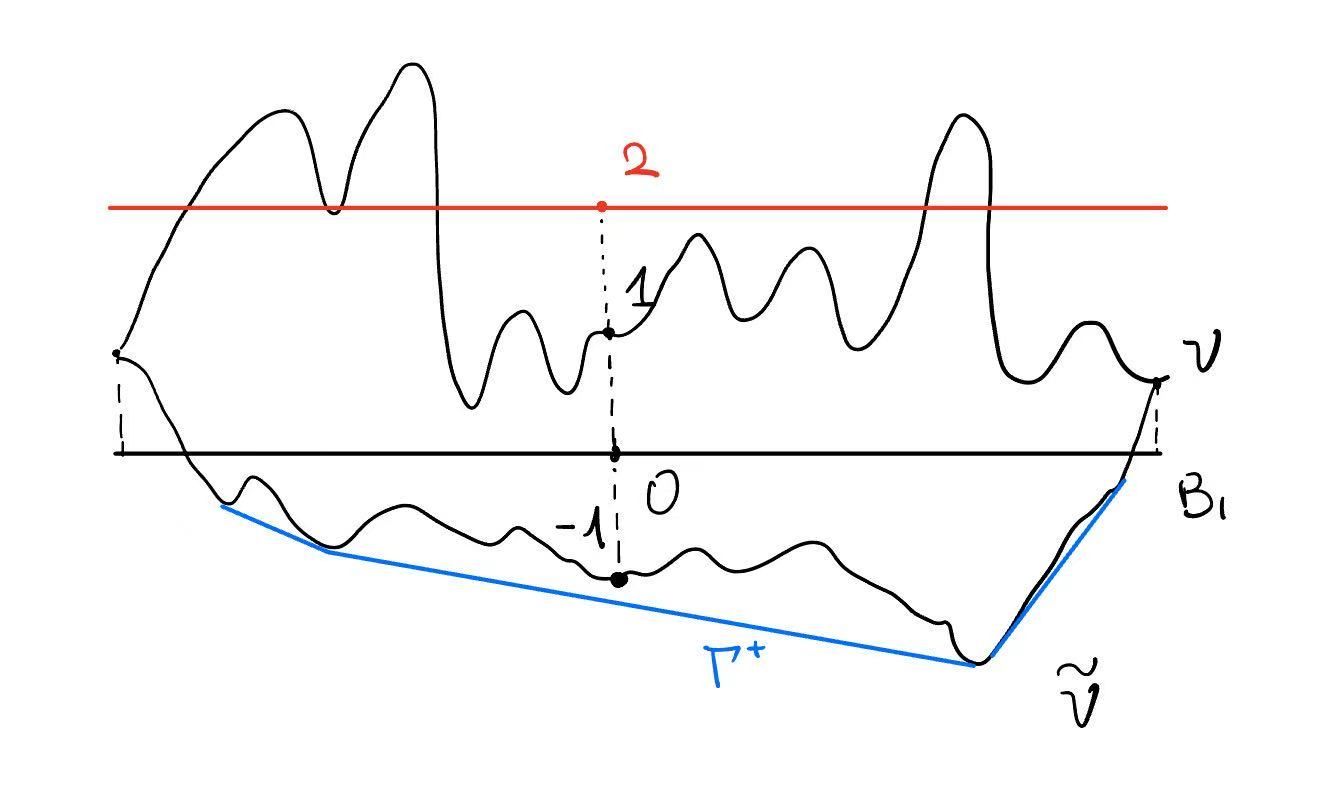}
\end{figure}

Applying the ABP estimate to $\tilde{v}$ yields
\[1\le \sup (\tilde{v}^-)^n\le \int_{\Gamma^+}C^n.\]
If $x_0\in \Gamma^+$, then $\tilde{v}(x_0)<0$ (for more preciser auxiliary functions, seeing the sublemma below), which implies $v(x_0)<2$. Hence $\Gamma^+\subset \{v<2\}$, and therefore
\[1\le C^n|\{v<2\}|.\]

We now state the general version of the measure estimate.

\begin{sublemma}[Contact estimate]\label{sublemma: measure estimate}
Let $x_0\in B_{1/4}$ and $M>0$ be a large constant. Then
\[|\{v<M\}\cap B_{1/8}(x_0)|\ge \theta|B_{1/8}(x_0)|,\]
where $\theta>0$ is universal.
\end{sublemma}

\begin{proof}
Construct an auxiliary function
\[
\gamma(x):=\begin{cases}
    \frac{\epsilon}{2}|x-x_0|^2-A,& \text{in } B_{1/8}(x_0),\\[4pt]
    -|x-x_0|^{-\sigma}-B,&\text{outside } B_{1/8}(x_0),
\end{cases}
\]
with constants $A,B$ chosen so that
\[
a_{ij}D_{ij}\gamma\le\begin{cases}
    C,& \text{in } B_{1/8}(x_0),\\[2pt]
    0,&\text{outside } B_{1/8}(x_0),
\end{cases}
\]
and set $\tilde{v}:=v+\gamma$ (see Figure~\ref{measure estimate 3}). We can arrange $\gamma$ so that $\tilde{v}(0)=-2$ and $\inf \gamma\ge -M$, without loss of generality, set $M>2$.

\begin{figure}[H]
    \centering
    \includegraphics[width=0.8\linewidth]{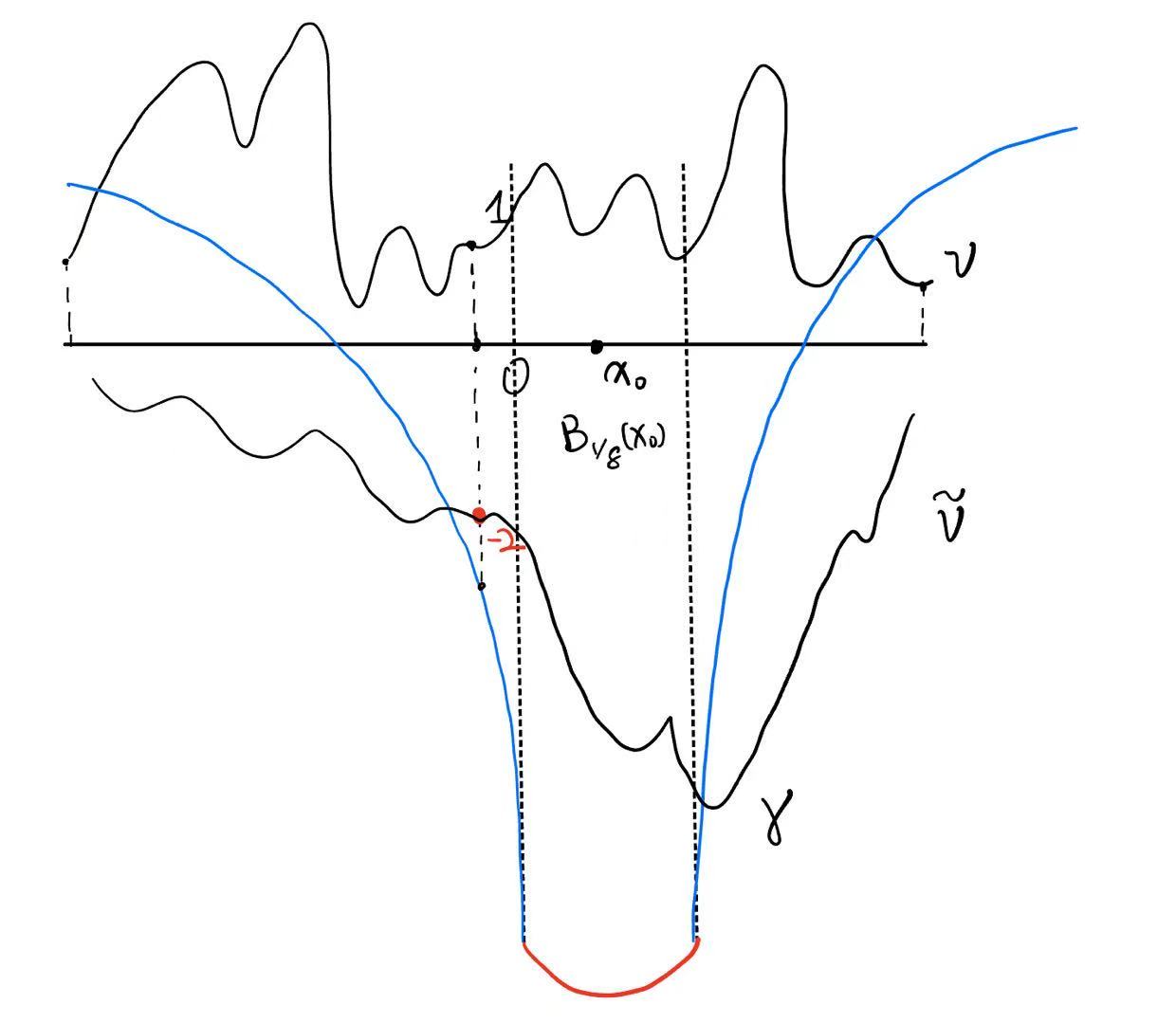}
    \caption{Measure estimate.}
    \label{measure estimate 3}
\end{figure}

Consequently,
\[
a_{ij}D_{ij}\tilde{v}\le\begin{cases}
    C,& \text{in } B_{1/8}(x_0),\\
    0,&\text{outside } B_{1/8}(x_0).
\end{cases}
\]
Applying the ABP estimate to $\tilde{v}$ yields
\[
2^n\le \sup (\tilde{v}^-)^n\le C^n|\Gamma^+|.
\]
If $x\in\Gamma^+$, then $\tilde{v}(x)<0$ (it is clear in $B_{1/8}(x_0)$; and $\Gamma^+$ can not escape from $B_{1/8}(x_0)$ because $\tilde{v}$ is a supersolution outside $B_{1/8}(x_0)$), which implies $v(x)<M$. Hence $\Gamma^+\subset\{v<M\}\cap B_{1/8}(x_0)$, and the desired inequality follows.
\end{proof}

The second step provides the weak $L^\epsilon$ estimate.

\begin{proof}[Proof of Lemma \ref{lemma: Krylov-Safonov weak Harnack inequality}]
Define
\[
A_k:=\{v>M^k\}\cap B_{1/2},\quad k=0,1,2,\cdots.
\]
We will show that the pair $(A_k, A_{k-1})$ satisfies the conditions of the A-B lemma with $\delta = 1-\theta$.  

Assume that for some cube $Q\subset B_{1/2}$ we have
\begin{equation}\label{inequality: measure argument 1}
\frac{|A_k\cap Q|}{|Q|} > \delta,
\end{equation}
which is equivalent to
\[
|\{v\le M^k\}\cap Q| < \theta|Q|.
\]
Let $\tilde{v}=v/M^{k-1}$ and rescale $Q$ to a cube $\mathbb{Q}_{1/8}(x_0)$. Then
\[
|\{\tilde{v}\le M\}\cap\mathbb{Q}_{1/8}(x_0)| < \theta|\mathbb{Q}_{1/8}(x_0)|.
\]
We claim that $\tilde{v}>1$ on the dilated cube $\widetilde{\mathbb{Q}}_{1/8}(x_0)$. If not, there exists $\tilde{x}\in\widetilde{\mathbb{Q}}_{1/8}(x_0)$ with $\tilde{v}(\tilde{x})\le1$. By Lemma~\ref{sublemma: measure estimate} applied to $\tilde{v}$ (after a translation), we obtain
\[
|\{\tilde{v}\le M\}\cap\mathbb{Q}_{1/8}(x_0)| \ge \theta|\mathbb{Q}_{1/8}(x_0)|,
\]
contradicting the previous inequality. Hence $\tilde{v}>1$ on $\widetilde{\mathbb{Q}}_{1/8}(x_0)$, which means $v>M^{k-1}$ there; thus $\widetilde{Q}\subset A_{k-1}$.

Therefore, the hypotheses of the A-B lemma are satisfied, and we obtain
\[
|A_k|\le (1-\theta)|A_{k-1}|.
\]
Iterating gives $|A_k|\le (1-\theta)^k|A_0|$. Setting $\epsilon = -\log_{M}(1-\theta)>0$, we have
\[
|\{v>M^k\}|\le M^{-\epsilon k},
\]
which implies $v\in L^\epsilon_{\text{weak}}(B_{1/2})$.
\end{proof}

\begin{remark}
One can explicitly estimate the "$L^\epsilon$-norm" in $B_{1/2}$:
\begin{align*}
    \int_{B_{1/2}} v^\epsilon &\le M^\epsilon|\{v<M\}|+M^{2\epsilon}|\{M\le v<M^2\}|+\cdots\\
    &\le M^\epsilon+M^{2\epsilon}(1-\theta)+\cdots+M^{k\epsilon}(1-\theta)^{k-1}+\cdots\\
    &\le \frac{M^\epsilon}{1-M^\epsilon(1-\theta)}<\infty.
\end{align*}
\end{remark}

\begin{corollary}[$L^\epsilon$-estimate]\label{corollary: Lepsilon estimate}
If $v\ge 0$ is a supersolution of \eqref{equation: fully nonlinear equation 4} in $\mathbb{Q}_3$, then
\[
\|v\|_{L^\epsilon(\mathbb{Q}_1)}\lesssim \inf_{\mathbb{Q}_2} v.
\]
\end{corollary}

\begin{remark}
Compare it with Lemma~\ref{lemma: Moser weak Harnack}.
\end{remark}

The above argument already yields $C^\alpha$ regularity. The proof of the next result follows the notes of Professor \href{https://sites.math.washington.edu/~yuan/}{Yu Yuan}.

\begin{corollary}[$C^\alpha$-estimate]\label{corollary: Calpha estimate}
Let $v$ be a nonnegative \textit{solution} of \eqref{equation: fully nonlinear equation 4}. Then
\[
\|v\|_{C^\alpha(B_{1/2})}\lesssim\|v\|_{L^\infty(B_1)}.
\]
\end{corollary}

\begin{proof}
It suffices to show
\[
\mathrm{osc}_{\mathbb{Q}_1} v \le \sigma\,\mathrm{osc}_{\mathbb{Q}_{2}} v,\quad \sigma<1.
\]
Define
\[
w:=\frac{v-\inf_{\mathbb{Q}_{2}} v}{\mathrm{osc}_{\mathbb{Q}_{2}} v},\quad 0\le w\le 1,
\]
which also solves \eqref{equation: fully nonlinear equation 4}. We consider two cases; if the first does not hold, we replace $w$ by $1-w$. Recall that $w$ is a \textit{solution}.
\begin{enumerate}
    \item[CASE I:] $|\{w\ge 1/2\}\cap \mathbb{Q}_1|\ge \dfrac{1}{2}$.
    \item[CASE II:] $|\{w\ge 1/2\}\cap \mathbb{Q}_1|< \dfrac{1}{2}$.
\end{enumerate}

Assume CASE I holds. Then
\[
\frac{1}{2}\left(\frac{|\mathbb{Q}_1|}{2}\right)^{1/\epsilon}\le 
\left(\int_{\mathbb{Q}_1} w^\epsilon\right)^{1/\epsilon}\lesssim 
\inf_{\mathbb{Q}_2} w \le \inf_{\mathbb{Q}_1} w.
\]
Hence $\inf_{\mathbb{Q}_1} w\ge c_0>0$.

\begin{figure}
    \centering
    \includegraphics[width=0.6\linewidth]{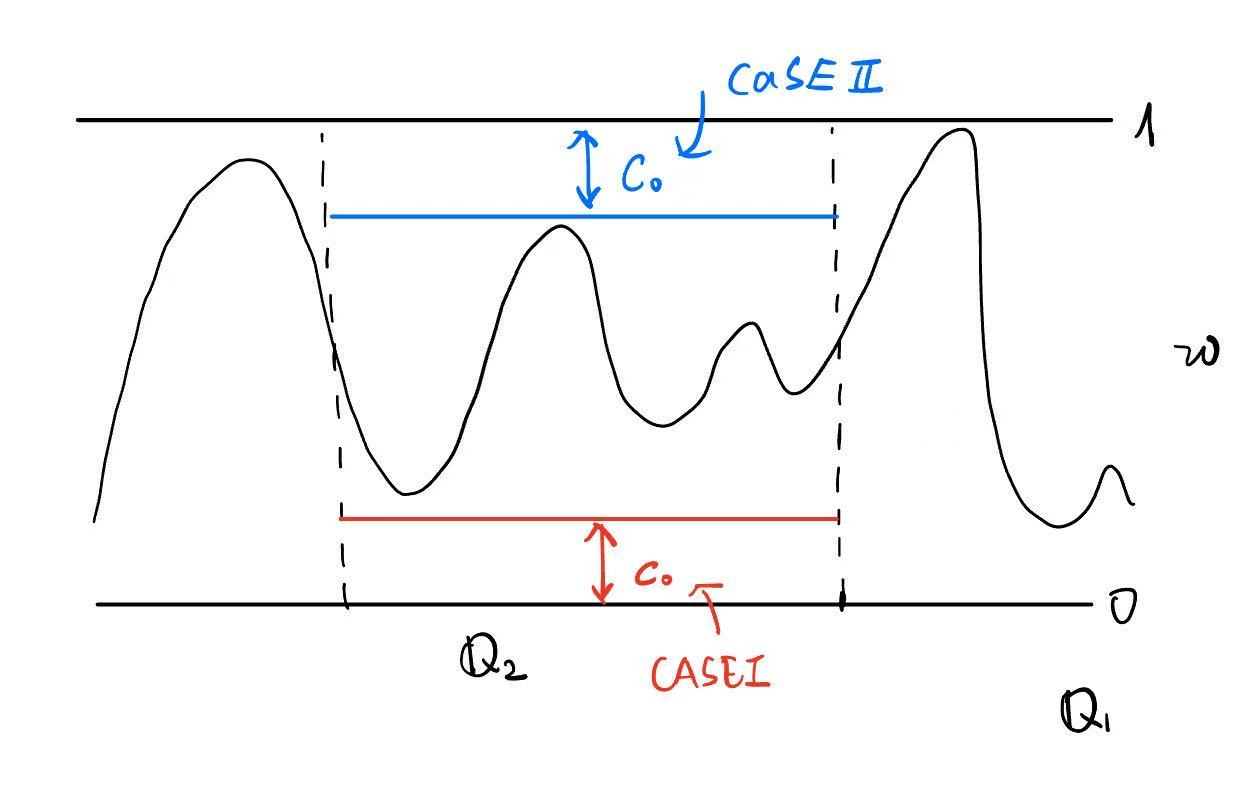}
    \caption{Oscillation decay}

\end{figure}
\end{proof}

\begin{remark}
Note that the proof fails if $v$ is merely a \textit{subsolution} or \textit{supersolution}; this constitutes the main difficulty in the Evans-Krylov theory. Also see Lemma~\ref{lemma: Moser weak Harnack}.
\end{remark}

\begin{remark}
The proof above naturally leads to a Liouville-type theorem analogous to Corollary~\ref{corollary: Liouville 1}.
\end{remark}

\subsection{Local maximum principle}
Finally, we present the second part of the proof, which implies Theorem \ref{theorem:KS-Harnack}.

\begin{lemma}[Local maximum principle]\label{lemma: local maximum principle}
Let $v$ be a \textit{subsolution} belonging to $L^\epsilon(B_1)$ with $\|v\|_{L^\epsilon(B_1)}\le1$. Then there exists a large constant $M>0$ such that
\[
\sup_{B_{1/2}} v \le M.
\]
\end{lemma}

\begin{remark}
Compare it with Lemma~\ref{lemma: DNM 1} and Lemma~\ref{lemma: moser iteration}.
\end{remark}

\begin{proof}
The proof relies on so-called \textit{blow‑up} argument. Assume, for contradiction, that there exists $x_0\in B_{1/2}$ with $v(x_0)=M_0$. We claim that there is a small cube $Q(x_0)\subset B_1$ such that
\[
\sup_{Q(x_0)} v \ge (1+\delta)M_0,\quad \text{for some universal $\delta$}.
\]

\begin{figure}[H]
    \centering
    \includegraphics[width=0.6\linewidth]{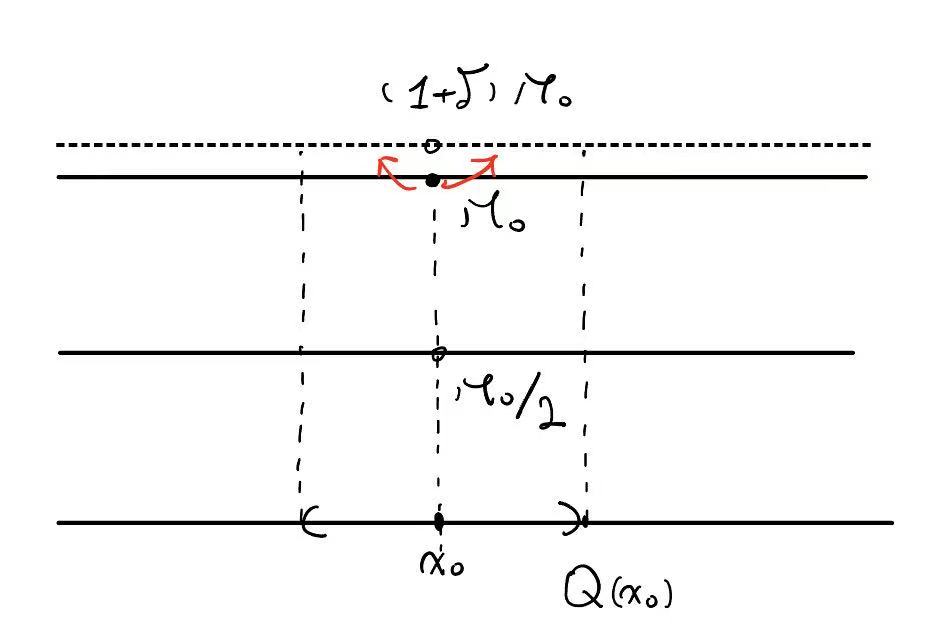}
    \caption{blow-up}
\end{figure}

Since $v\in L^\epsilon(B_1)$ and $\|v\|_{L^\epsilon}\le1$,
\[
\left|\left\{v>\frac{M_0}{2}\right\}\cap \mathbb{Q}_1\right| \le \left(\frac{M_0}{2}\right)^{-\epsilon}.
\]
Choose $Q(x_0)$ with $|Q(x_0)|>2(\frac{M_0}{2})^{-\epsilon}$ (recall $M$ may be taken sufficiently large); then
\[
\left|\left\{v>\frac{M_0}{2}\right\}\cap Q(x_0)\right| \le \frac12|Q(x_0)|.
\]

Now prove the claim. If not, i.e., $v< (1+\delta)M_0$ in $Q(x_0)$. Define
\[
w:=\frac{(1+\delta)M_0-v}{\delta M_0},
\]
which is a nonnegative supersolution with $w(x_0)=1$. By K-S weak Harnack inequality (Lemma~\ref{lemma: Krylov-Safonov weak Harnack inequality}), $w\in L^\epsilon_{weak}(Q(x_0))$. Note that, if $v\le M_0/2$, then
\[
w \ge \frac{(1+\delta)M_0-M_0/2}{\delta M_0}= \frac{\frac12+\delta}{\delta}.
\]
Hence
\begin{align*}
\left|\left\{v\le\frac{M_0}{2}\right\}\cap Q(x_0)\right|
&\le \left|\left\{w\ge\frac{\frac12+\delta}{\delta}\right\}\cap Q(x_0)\right| \\
&\le \left(\frac{\delta}{\frac12+\delta}\right)^{\!\epsilon}|Q(x_0)| < \frac12|Q(x_0)|,
\end{align*}
for sufficiently small $\delta$. This contradiction establishes the claim.

Hence, one may choose a $x_1\in Q(x_0)$ such that
\[v(x_1)\ge (1+\delta)M_0,\quad |x_1-x_0|\lesssim M_0^{-\epsilon};\]
\[v(x_2)\ge (1+\delta)^2M_0,\quad |x_2-x_1|\lesssim [(1+\delta)M_0]^{-\epsilon};\]
\[v(x_k)\ge(1+\delta)^kM_0,\quad |x_k-x_{k-1}|\lesssim [(1+\delta)^kM_0]^{-\epsilon},\quad k=3,4,\cdots.\]
Note that
\begin{align*}
    |x_k-x_0|&\le |x_k-x_{k-1}|+\cdots+|x_1-x_0|\\
    &\le \left[1+(1+\delta)+\cdots+(1+\delta)^k\right]^{-\epsilon}M_0^{-\epsilon}\\
    &\le \left(\frac{\delta}{(1+\delta)^k-1}\right)^\epsilon M_0^{-\epsilon}<1/8,
\end{align*}
as $\delta$ sufficiently small; and $v(x_k)\to \infty$ as $k\to\infty$, which leads to the contradiction.
\end{proof}
\begin{remark}
This yields a contradiction because, starting from a point $x_k$ with large value, one can iteratively produce points arbitrarily close to $x_0$ with ever larger values, contradicting the definition of $L^\epsilon_{weak}$.
\end{remark}

\section{Concave property and Evans-Krylov  theory}
\subsection{Concave assumption of $F$ and Gap}
We have established the Krylov-Safonov theory, i.e., the $C^{1,\alpha}$ estimate for
\[F_{ij}(D^2u)D_{ij}u_e=0.\]
To obtain $C^{2,\alpha}$ estimates, an additional assumption is required, namely that $F$ is \textit{concave}. Differentiating the equation again in a direction $f\in\mathbb{S}^{n-1}$ yields, formally,
\[F_{ij,kl}(D^2u)D_{ij}u_eD_{kl}u_f+F_{ij}(D^2u)D_{ij}u_{ef}=0.\]
Since $F$ is concave, we have
\begin{equation}\label{equation: fully nonlinear equation 5}
    F_{ij}(D^2u)D_{ij}v\le 0,\quad v:=u_{ef}.
\end{equation}
Thus, $v$ is a \textit{subsolution}. Recall the $L^\epsilon$ estimate in Krylov-Safonov theory (Corollary~\ref{corollary: Lepsilon estimate}); however, there is a gap because we lack a corresponding supersolution. The fact that $u_{ef}\in L^\epsilon$ was proved separately by Evans \cite{evans1985some} and Lin \cite{lin1986second}. Then, by the local maximum principle (Lemma~\ref{lemma: local maximum principle}), we obtain $u\in C^{1,1}$. For further details, we refer the reader to \cite{caffarelli1995fully}.

It remains to prove that $u\in C^{2,\alpha}$, i.e., $D^2u\in C^\alpha$. The $C^\alpha$ estimate (Corollary~\ref{corollary: Calpha estimate}) requires a "solution", but we only know that $v$ is a subsolution. This gap constitutes the main difficulty in the Evans-Krylov theory. Fortunately, however, the relation
\[(D^2u)^+\sim(D^2u)^-\Longleftrightarrow v^+\sim v^-\]
suggests that $v$ is almost a \textit{supersolution}, an idea that will be elaborated in the following arguments.

Firstly, we state the main result.
\begin{theorem}[Evans-Krylov, 1982 \cite{evans1982classical,krylov1983boundedly,krylov1984boundedly}]
Let $u$ be a $C^{1,1}$ solution to
\[F(D^2u)=0,\quad \text{in}\ B_1,\]
where $F$ is uniformly elliptic and \textit{concave} (or convex). Then $u\in C^{2,\alpha}(B_{1/2})$, precisely
\[\|u\|_{C^{2,\alpha}(B_{1/2})}\lesssim \|u\|_{L^\infty(B_1)}.\]
\end{theorem}

Set $\mathcal{H}:=D^2u:B_1\to \mathcal{S}$ as the Hessian map, where $\mathcal{S}$ denotes the space of real symmetric matrices with the usual Euclidean norm. The theorem follows from the following oscillation lemma. Without loss of generality, assume $\mathrm{diam}\ \mathcal{H}(B_1)=1$.
\begin{lemma}[Oscillation decay]\label{lemma: oscillation decay 2}
There exists a $\rho_0<1$ such that
\[\mathrm{diam}\ \mathcal{H}(B_{\rho_0})<\frac{1}{2}.\]
\end{lemma}

\subsection{Deletion argument}
The idea relies on a \textit{deletion} argument: one can remove a small ball from the covering of $\mathcal{H}(B_1)$ without affecting the covering of $\mathcal{H}(B_{1/2})$. This is made precise by the covering lemma stated below.

Let \(\mathbf{t}:=\sup_{B_1}\lambda_{\max}(D^{2}u)\).

\begin{sublemma}[Covering]\label{sublemma:covering}
There exists a universal \(0<\delta\ll1\) such that \(\mathcal{H}(B_{1})\) can be covered by finitely many balls \(\{B_{\delta}(M_{i})\}\). Moreover, for some direction \(g\in\mathbb{S}^{n-1}\) and some ball \(B_{\delta}(M_{3})\) we have
\[
\mathbf{t}-u_{gg}(x)\gtrsim\frac{1}{8},\quad\forall x\in\mathcal{H}^{-1}\left(B_{\delta}(M_{3})\right),
\]
and
\[
\left|\mathcal{H}^{-1}\left(B_{\delta}(M_{3})\right)\right|\gtrsim\delta^{\,n^{2}}.
\]
\end{sublemma}

\begin{proof}
Cover \(\mathcal{H}(B_{1})\) by balls of radius \(\delta\), which will be determined at the end of proof; their number satisfies \(\#\lesssim\delta^{-n^{2}}\). Hence some ball \(B_{\delta}(M_{3})\) satisfies the measure estimate
\begin{equation}\label{equation: measure estimate 4}
\left|\mathcal{H}^{-1}\left(B_{\delta}(M_{3})\right)\right|\gtrsim\delta^{\,n^{2}}.
\end{equation}

Since \(\operatorname{diam}\mathcal{H}(B_{1})=1\), there exists \(M\in\mathcal{H}(B_{1})\) with \(\|M-M_{3}\|\ge\frac14\). Both \(F(M)=F(M_{3})=0\). By uniform ellipticity,
\[
\lambda\|(M-M_{3})^{+}\|-\Lambda\|(M-M_{3})^{-}\|\le0\le\Lambda\|(M-M_{3})^{+}\|-\lambda\|(M-M_{3})^{-}\|,
\]
which implies \(\|(M-M_{3})^{+}\|\sim\|(M-M_{3})^{-}\|\). Consequently,
\[
\|(M-M_{3})^{+}\|\gtrsim\frac18,\quad\|(M-M_{3})^{-}\|\gtrsim\frac18.
\]

Take points \(x_{3},y\in B_{1}\) with \(D^{2}u(x_{3})=M_{3}\) and \(D^{2}u(y)=M\). Choose a unit vector \(g\) (e.g., an eigenvector of \((M-M_{3})^{+}\) corresponding to its largest eigenvalue) such that
\[
u_{gg}(y)-u_{gg}(x_{3})\gtrsim\frac18.
\]
Since \(\mathbf{t}\ge u_{gg}(y)\),
\[
\mathbf{t}\gtrsim u_{gg}(x_{3})+\frac18.
\]

For any \(z\in\mathcal{H}^{-1}(B_{\delta}(M_{3}))\), we have \(|u_{gg}(z)-u_{gg}(x_{3})|\lesssim \delta\). Choosing \(\delta\) sufficiently small (universally) gives
\[
\mathbf{t}-u_{gg}(z)\gtrsim\frac18-\delta\gtrsim\frac18,
\]
while preserving \eqref{equation: measure estimate 4}.
\end{proof}

\begin{figure}[H]
    \centering
    \includegraphics[width=0.9\linewidth]{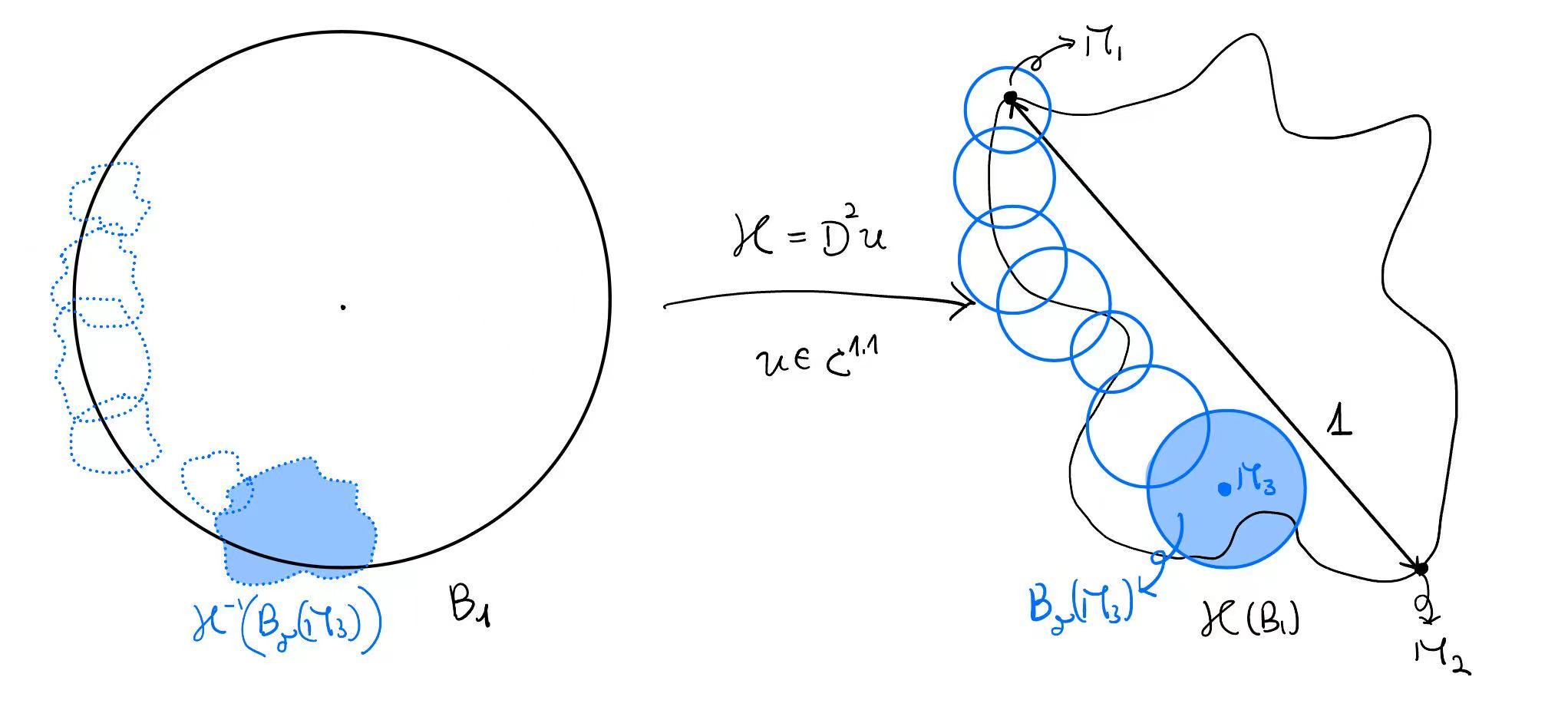}
\end{figure}

Now consider the function
\[
\tilde{v}:=\mathbf{t}-u_{gg},
\]
which satisfies:
\begin{enumerate}
    \item[(i)] $\tilde{v}\ge0$ in $B_1$;
    \item[(ii)] $\tilde{v}$ is a \textit{supersolution};
    \item[(iii)] $\tilde{v}\gtrsim\frac18$ on $\mathcal{H}^{-1}\left(B_\delta(M_3)\right)$.
\end{enumerate}
By (i), (ii) and the $L^\epsilon$ estimate (Corollary~\ref{corollary: Lepsilon estimate}),
\[
\inf_{B_{1/2}}\tilde{v}\gtrsim\left(\int_{B_1}\tilde{v}^\epsilon\right)^{1/\epsilon}
\gtrsim\frac18\,\delta^{\,n^{2}/\epsilon}=:\theta,
\]
where the last bound uses (iii). Hence
\[
u_{gg}(x)\le\mathbf{t}-\theta\quad\forall x\in B_{1/2}.
\]

We now prove Lemma~\ref{lemma: oscillation decay 2} via a \emph{deletion} argument.

\begin{proof}[Proof of Lemma~\ref{lemma: oscillation decay 2}]
First, cover $\mathcal{H}(B_1)$ by finitely many balls of radius $\theta$ (which is smaller than the earlier $\delta$). The estimate above shows that $u_{gg}(x)\le\mathbf{t}-\theta$ on $B_{1/2}$; therefore $\mathcal{H}(B_{1/2})$ can be covered by the same family \emph{without} the ball $B_\theta(\mathbf{t})$ (see Figure~\ref{deletion_1}). Thus one ball is effectively deleted from the covering.

\begin{figure}[H]
    \centering
    \includegraphics[width=0.9\linewidth]{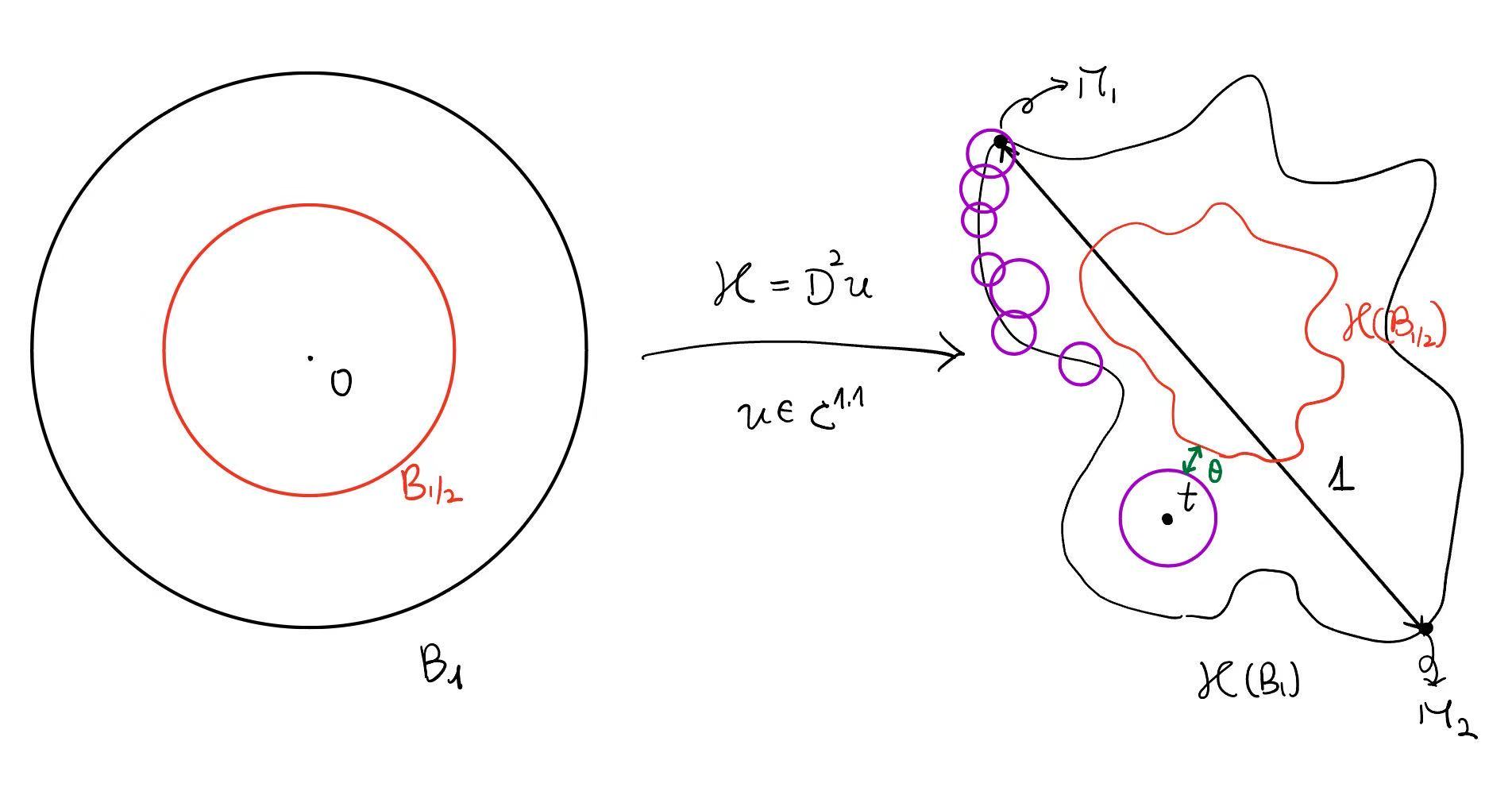}
    \caption{deletion}
    \label{deletion_1}
\end{figure}

If after this deletion we have $\operatorname{diam}\mathcal{H}(B_{1/2})\le\frac12$, we may take $\rho_0=1/2$ and the lemma is proved. Otherwise,
$\operatorname{diam}\mathcal{H}(B_{1/2})>\frac12$, and we can apply Sublemma~\ref{sublemma:covering} to the rescaled function $v_1(x):=2^2\,v(x/2)$. Since the number of balls in the covering is finite, after iterating this procedure a finite number $k_0$ of times we obtain
\[
\operatorname{diam}\mathcal{H}(B_{2^{-k_0}})<\frac12.
\]
Taking $\rho_0=2^{-k_0}$ completes the proof.
\end{proof}

\subsection{Notes}
\subsubsection{Krylov-Safonov}
A simpler proof of the Krylov–Safonov theorem is provided in \cite{mooney2019proof}, whose approach does not rely on the \textit{localization} technique.

\subsubsection{Evans–Krylov}
Since the Evans–Krylov theory requires concavity of $F$, a natural question is whether solutions are always $C^2$. The answer is affirmative for $n=2$ \cite{nirenberg1953nonlinear} (see also Theorem 4.8 of \cite{fernandez2023regularity}), where solutions are in fact $C^{2,\alpha}$. Counterexamples for $n\ge5$ were constructed by Nadirashvili-Vladuts \cite{nadirashvili2007nonclassical,nadirashvili2008singular,nadirashvili2013singular}. The cases $n=3,4$ remain open. We refer to \cite{caffarelli1995fully} for a comprehensive treatment of fully nonlinear equations. Recently, Caffarelli-Yuan \cite{caffarelli2000priori} established $C^{2,\alpha}$ regularity under a convexity assumption on the level set $\{F(D^2u)=0\}$; Cabré-Caffarelli \cite{cabre2003interior} obtained the same results for certain nonconvex operators $F$. Later, Collins \cite{collins2016c} and Goffi \cite{goffi2024high} studied regularity under weaker convexity conditions.

\subsubsection{Non-uniformly elliptic equations}
The equations discussed in this note are all \textit{uniformly} elliptic. In fact, many classical equations are \textit{locally} uniformly elliptic but not \textit{globally} uniformly elliptic, such as the Monge–Ampère equation $\mathrm{det}D^2u=0$ (or, more generally, the $\sigma_k$-Hessian equations), and the minimal surface equation $\mathrm{div}\left(\frac{\nabla u}{\sqrt{1+|\nabla u|^2}}\right)=0$ (or, more generally, the $\sigma_k$-curvature equations). The main difficulty in studying locally uniformly elliptic equations is that \textit{scaling} can alter the uniform ellipticity. Consequently, the \textit{renormalization} and \textit{blow-up} arguments commonly used in previous proofs are no longer applicable.

Fortunately, to prove the weak Harnack inequality—for example, see Lemma~\ref{lemma: Krylov-Safonov weak Harnack inequality}—the so-called \textit{sliding paraboloids} method, which estimates the contact measure by approximating a series of paraboloids, avoids the self-rescaling that would disturb uniform ellipticity. This method was first introduced by Cabré \cite{cabre1997nondivergent}, who extended the Krylov–Safonov estimate to manifolds with negative sectional curvature. Later, Savin \cite{savin2009regularity}, inspired also by Caffarelli-Córdoba \cite{caffarelli1993elementary}, used this approach to improve De Giorgi’s original flatness theorem \cite{giusti1984minimal} (and, more generally, Allard regularity \cite{allard1972first}), thereby proving the De Giorgi conjecture under a technical condition in dimensions $n\le 8$. Furthermore, Savin \cite{savin2007small} extended the idea to locally uniformly elliptic equations with \textit{small perturbations} (or almost flatness), building on work of Caffarelli \cite{caffarelli1987harnack} and Caffarelli–Wang \cite{caffarelli1993harnack} (see further details in \cite{caffarelli1995fully}). Subsequently, dos Prazeres and Teixeira \cite{dos2016asymptotics} treated the \textit{non-homogeneous} case under uniform ellipticity, while Fan \cite{fan2025generalization} generalized it to the non-homogeneous case with only local uniform ellipticity. Wang \cite{wang2013small} studied the parabolic case, and Yu \cite{yu2017small} investigated the nonlocal case.

\bibliographystyle{alpha}
\bibliography{arxiv}
\end{document}